\documentclass{article}
\usepackage{amssymb}
\usepackage{amsfonts}
\usepackage{amsmath}

\newtheorem{theorem}{Theorem}

\newtheorem{corollary}[theorem]{Corollary}

\newtheorem{definition}[theorem]{Definition}

\newtheorem{lemma}[theorem]{Lemma}

\newtheorem{proposition}[theorem]{Proposition}
\newtheorem{remark}[theorem]{Remark}

\newenvironment{proof}[1][Proof]{\noindent\textbf{#1.} }{\ \rule{0.5em}{0.5em}}
\input{tcilatex}

\begin{document}

\title{Affine $ADE$ bundles over \\
complex surfaces with $\mathit{p}_{g}=0$}
\author{Yunxia Chen \& Naichung Conan Leung}
\date{}
\maketitle

\begin{abstract}
We study simply-laced simple affine Lie algebra bundles over complex
surfaces $X$. Given any Kodaira curve $C$ in $X$, we construct such a bundle
over $X$. After deformations, it becomes trivial on every irreducible
component of $C$ provided that $\mathit{p}_{g}(X)=0$.

When $X$ is a blowup of $\mathbb{P}^{2}$ at nine points, there is a
canonical $\hat{E}_{8}$-bundle $\mathcal{E}_{0}^{\hat{E}_{8}}$ over $X$. We
show that the geometry of $X$ can be reflected by the deformability of $%
\mathcal{E}_{0}^{\hat{E}_{8}}$.
\end{abstract}

\section{Introduction}

Given a complex surface $X$ and a sublattice $\Lambda \subset Pic\left(
X\right) $, if $\Lambda $ is isomorphic to the root lattice $\Lambda _{%
\mathfrak{g}}$ of a simple Lie algebra $\mathfrak{g}$, then we have a root
system $\Phi $ of $\mathfrak{g}$ and we can associate a Lie algebra bundle $%
\mathcal{E}_{0}^{\mathfrak{g}}$ over $X$ \cite{FM}\cite{LZ}\cite{LZ2},%
\begin{equation*}
\mathcal{E}_{0}^{\mathfrak{g}}:=O_{X}^{\oplus r}\oplus \bigoplus_{\alpha \in
\Phi }O_{X}(\alpha )\text{.}
\end{equation*}%
This can be generalized to affine Lie algebras $\widehat{\mathfrak{g}}$ \cite%
{LXZ}.

There are many instances that this happens: when $X_{n}$ is a del Pezzo
surface, namely a blowup of $\mathbb{P}^{2}$ at $n\leq 8$ points in general
position (or $\mathbb{P}^{1}\times \mathbb{P}^{1}$), $\left\langle
K_{X_{n}}\right\rangle ^{\bot }\subset Pic\left( X_{n}\right) $ is
isomorphic to $\Lambda _{E_{n}}$. Thus we have an $E_{n}$-bundle over $X_{n}$%
. By restriction, we have an $E_{n}$-bundle over any anti-canonical curve $%
\Sigma $ in $X_{n}$. Notice that $\Sigma $ is always an elliptic curve. For
a fixed elliptic curve $\Sigma $, the above construction gives a bijection
between del Pezzo surfaces containing $\Sigma $ and $E_{n}$-bundles over $%
\Sigma $ \cite{D}\cite{D2}\cite{FMW}\cite{LZ}\cite{L}\cite{XZ}. Such an
identification was predicted by the F-theory/string duality in physics \cite%
{FMW}. This was generalized to all simple Lie algebras in \cite{LZ}\cite{LZ2}%
. When $n=9$, $X_{9}$ is not Fano and $E_{9}=\hat{E}_{8}$ is an affine Lie
algebra. Corresponding results for the $\hat{E}_{8}$-bundle over $X_{9}$ are
obtained in \cite{LXZ}.

When $X$ is the canonical resolution of a surface $X^{\prime }$ with a
rational singularity, the exceptional curve $C=\cup C_{i}$ is an $ADE$ curve
of type $\mathfrak{g}$. Therefore $C_{i}$'s span a sublattice of $Pic\left(
X\right) $ which is isomorphic to $\Lambda _{\mathfrak{g}}$, thus giving a $%
\mathfrak{g}$-bundle $\mathcal{E}_{0}^{\mathfrak{g}}$ over $X$. When $%
p_{g}\left( X\right) =0$, there exists a deformation $\mathcal{E}_{\varphi
}^{\mathfrak{g}}$ of $\mathcal{E}_{0}^{\mathfrak{g}}$ such that $\mathcal{E}%
_{\varphi }^{\mathfrak{g}}$ is trivial on each $C_{i}$, thus it descends to
the singular surface $X^{\prime }$ \cite{CL}\cite{FM}.

When $X$ is a relatively minimal elliptic surface, Kodaira classified all
possible singular fibers (see e.g. \cite{BPV}) and we call such a curve $%
C=\cup C_{i}$ a Kodaira curve. Its irreducible components $C_{i}$'s span a
sublattice of $Pic\left( X\right) $ which is isomorphic to the root lattice
of an \textit{affine} root system $\Phi _{\widehat{\mathfrak{g}}}$ and
therefore we can construct an affine Lie algebra bundle $\mathcal{E}_{0}^{%
\widehat{\mathfrak{g}}}$ over $X$.

\begin{theorem}
$($Lemma \ref{cohomology}, \textit{Proposition \ref{holomorphic2}} and
Theorem \ref{thm2}$)$ Given any complex surface $X$ with $p_{g}=0$. If $X$
has a Kodaira curve $C=\cup _{i=0}^{r}C_{i}$ of type $\widehat{\mathfrak{g}}$%
, then

$(i)$ given any $(\varphi _{C_{i}})_{i=0}^{r}\in \Omega
^{0,1}(X,\bigoplus_{i=0}^{r}O(C_{i}))$ with $\overline{\partial }\varphi
_{C_{i}}=0$ for every $i$, it can be extended to $\varphi =(\varphi _{\alpha
})_{\alpha \in \Phi _{\widehat{\mathfrak{g}}}^{+}}\in \Omega
^{0,1}(X,\bigoplus_{\alpha \in \Phi _{\widehat{\mathfrak{g}}}^{+}}O(\alpha
)) $ such that $\overline{\partial }_{\varphi }:=\overline{\partial }%
+ad(\varphi )$ is a holomorphic structure on $\mathcal{E}_{0}^{\widehat{%
\mathfrak{g}}}$. We denote the new bundle as $\mathcal{E}_{\varphi }^{%
\widehat{\mathfrak{g}}}$.

$(ii)$ $\overline{\partial }_{\varphi }$ is compatible with the Lie algebra
structure on $\mathcal{E}_{0}^{\widehat{\mathfrak{g}}}$.

$(iii)$ $\mathcal{E}_{\varphi }^{\widehat{\mathfrak{g}}}$ is trivial on $%
C_{i}$ if and only if $[\varphi _{C_{i}}|_{C_{i}}]\neq 0\in
H^{1}(C_{i},O_{C_{i}}(C_{i}))\cong \mathbb{C}$.

$(iv)$ There exists $[\varphi _{C_{i}}]\in H^{1}(X,O(C_{i}))$ such that $%
[\varphi _{C_{i}}|_{C_{i}}]\neq 0$.
\end{theorem}

In the second half of this paper, we explain how the geometry of $X_{9}$, a
blowup of $\mathbb{P}^{2}$ at nine points, can be reflected by the
deformability of the $\widehat{E}_{8}$-bundle $\mathcal{E}_{0}^{\hat{E}_{8}}$
over it. Among other things, we obtained the following results.

\begin{theorem}
$($Theorem \ref{general position}$)$ $\mathcal{E}_{0}^{\widehat{E}_{8}}$ is
totally non-deformable if and only if the nine blowup points in $\mathbb{P}%
^{2}$ are in general position.
\end{theorem}

\begin{theorem}
$($Theorem \ref{elliptic}$)$ Suppose $-K_{X_{9}}$ is nef, then

$(i)$ $X_{9}$ admits an elliptic fibration with a multiple fiber of
multiplicity $m$ $\left( m\geq 1\right) $ if and only if $\mathcal{E}_{0}^{%
\widehat{E}_{8}}$ is deformable in $(-mK)$-direction but not in $(-m+1)K$%
-direction.

$(ii)$ $X_{9}$ has a $($maximal$)$ $ADE$ curve $C$ of type $\mathfrak{g}$ if
and only if $\mathcal{E}_{0}^{\widehat{E}_{8}}$ is $($maximal$)$ $\mathfrak{g%
}$-deformable.

$(iii)$ $X_{9}$ has a $($maximal$)$ Kodaira curve $C$ of type $\widehat{%
\mathfrak{g}}$ if and only if $\mathcal{E}_{0}^{\widehat{E}_{8}}$ is $($%
maximal$)$ $\widehat{\mathfrak{g}}$-deformable.
\end{theorem}

The organization of this paper is as follows. Section $2$ gives the
construction of the (affine) $ADE$ Lie algebra bundles directly from
(affine) $ADE$ curves. In section $3$, we assume $\mathit{p}_{g}(X)=0$. We
construct deformations of the holomorphic structures on these bundles such
that the new bundles are trivial over irreducible components of the curve.
We will consider the $E_{n}$-bundle over a blowup of $\mathbb{P}^{2}$ at $%
n\leq 9$ points in section $4$ and show how the deformability of this bundle
can reflect the geometry of the underlying surface. In the appendix, we
review the basic construction of affine Lie algebras.

Notations: for a holomorphic bundle $(E_{0},\overline{\partial }_{0})$ with $%
E_{0}=\oplus _{i}O(D_{i})$, $\overline{\partial }_{0}$ means the $\overline{%
\partial }$-operator for the direct sum holomorphic structure. If we
construct a new holomorphic structure $\overline{\partial }_{\varphi }$ on $%
E_{0}$, we denote the resulting bundle as $E_{\varphi }$.

\bigskip

$\mathbf{Acknowledgements.}$ We are grateful to R.Friedman, E. Looijenga and
J.J. Zhang for many useful comments and discussions. The work of the second
author was supported by a research grant from the Research Grants Council of
the Hong Kong Special Administrative Region, China (reference No. 401411).

\section{Affine $ADE$ bundles from affine $ADE$ curves}

\subsection{$ADE$ and affine $ADE$ curves}

\begin{definition}
A curve $C=\cup C_{i}$ in a surface $X$ is called an $ADE$ $($resp. affine $%
ADE)$ curve of type $\mathfrak{g}$ $($resp. $\widehat{\mathfrak{g}})$ if
each $C_{i}$ is a smooth $(-2)$-curve in $X$ and the dual graph of $C$ is a
Dynkin diagram of the corresponding type.
\end{definition}

It is known that $C$ is an $ADE$ curve if and only if $C$ can be contracted
to a rational singularity. In this case, the intersection matrix $%
(C_{i}\cdot C_{j})<0$ \cite{BPV}.

If $C$ is an affine $ADE$ curve, then the intersection matrix $(C_{i}\cdot
C_{j})\leq 0$ and there exists unique $n_{i}$'s up to overall scalings such
that $F:=\sum n_{i}C_{i}$ satisfies $F\cdot F=0$. Dynkin diagrams of affine $%
ADE$ types are drawn as follows and the corresponding $n_{i}C_{i}$'s are
labelled in the pictures. $ADE$ Dynkin diagrams can be obtained by removing
the node corresponding to $C_{0}$.

\bigskip 
\begin{equation}
\setlength{\unitlength}{0.9cm}\begin{picture}(5, 2) \put(-2,
1){$\hat{A}_n:$}\put(0, 1){\circle*{.2}} \put(1.2, 1){\circle*{.2}}
\put(2.4, 1){\circle*{.2}} \put(3.6, 1){\circle*{.2}} \put(4.8,
1){\circle*{.2}} \put(2.4, 2.2){\circle*{.2}} \put(0,1){\line(2, 1){2.4}}
\put(4.8,1){\line(-2, 1){2.4}} \put(0,1){\line(1, 0){1.2}} \put(1.5,
1){\circle*{.1}} \put(1.8, 1){\circle*{.1}}\put(2.1, 1){\circle*{.1}}
\put(2.4, 1){\line(1, 0){1.2}} \put(3.6, 1){\line(1, 0){1.2}} \put(-0.5,
.5){$1 C_1$} \put(1, .5){$1 C_2$} \put(2.3, .5){$1 C_{n-2}$} \put(3.5,
.5){$1 C_{n-1}$} \put(4.6, .5){$1 C_{n}$} \put(2.6,2.1){$1 C_{0}$}
\end{picture}  \notag
\end{equation}%
\begin{equation*}
\setlength{\unitlength}{0.9cm}\begin{picture}(5, 3) \put(-2,
1){$\hat{D}_n:$} \put(0, 1){\circle*{.2}} \put(1.2, 1){\circle*{.2}}
\put(2.4,1){\circle*{.2}} \put(3.6,1){\circle*{.2}}
\put(4.8,1){\circle*{.2}} \put(3.6,2.2){\circle*{.2}}
\put(1.2,2.2){\circle*{.2}} \put(1.2, 1){\line(0, 1){1.2}} \put(1.5,
1){\circle*{.1}} \put(1.8, 1){\circle*{.1}} \put(2.1, 1){\circle*{.1}}
\put(0, 1){\line(1, 0){1.2}} \put(2.4, 1){\line(1, 0){1.2}} \put(3.6,
1){\line(1, 0){1.2}} \put(3.6, 1){\line(0, 1){1.2}} \put(-0.5, .5){$1 C_1$}
\put(1, .5){$2 C_2$} \put(2.3, .5){$2 C_{n-3}$} \put(3.5, .5){$2 C_{n-2}$}
\put(4.5, .5){$1 C_{n-1}$} \put(3.8, 2.1){$1 C_{n}$} \put(1.4, 2.1){$1
C_{0}$} \end{picture}
\end{equation*}%
\begin{equation*}
\setlength{\unitlength}{0.8cm}\begin{picture}(6, 3) \put(-2,
1){$\hat{E}_6:$}\put(0, 1){\circle*{.2}} \put(1.2, 1){\circle*{.2}}
\put(2.4, 1){\circle*{.2}} \put(3.6, 1){\circle*{.2}} \put(4.8,
1){\circle*{.2}} \put(2.4,2.0){\circle*{.2}} \put(2.4,3.0){\circle*{.2}}
\put(0, 1){\line(1, 0){1.2}} \put(1.2, 1){\line(1, 0){1.2}} \put(2.4,
1){\line(1, 0){1.2}} \put(3.6, 1){\line(1, 0){1.2}} \put(2.4, 1){\line(0,
1){1.0}} \put(2.4,2.0){\line(0, 1){1.0}} \put(-0.5, .5){$1 C_1$} \put(1.1,
.5){$2 C_2$} \put(2.3, .5){$3 C_{3}$} \put(3.5, .5){$2 C_{4}$} \put(4.7,
.5){$1 C_{5}$} \put(2.6, 1.9){$2 C_{6}$} \put(2.6, 2.9){$1 C_{0}$}
\end{picture}
\end{equation*}%
\begin{equation*}
\setlength{\unitlength}{0.8cm}\begin{picture}(6, 3) \put(-2,
1){$\hat{E}_7:$}\put(0, 1){\circle*{.2}} \put(1.2, 1){\circle*{.2}}
\put(2.4, 1){\circle*{.2}} \put(3.6, 1){\circle*{.2}} \put(4.8,
1){\circle*{.2}} \put(6, 1){\circle*{.2}} \put(7.2,1){\circle*{.2}}
\put(3.6,2.2){\circle*{.2}} \put(0, 1){\line(1, 0){1.2}} \put(1.2,
1){\line(1, 0){1.2}} \put(2.4, 1){\line(1, 0){1.2}} \put(3.6, 1){\line(1,
0){1.2}} \put(4.8,1){\line(1,0){1.2}} \put(6,1){\line(1,0){1.2}} \put(3.6,
1){\line(0, 1){1.2}} \put(-0.5, .5){$1 C_1$} \put(1.1, .5){$2 C_2$}
\put(2.3, .5){$3 C_{3}$} \put(3.5, .5){$4 C_{4}$} \put(4.7, .5){$3 C_{5}$}
\put(5.8, .5){$2 C_{6}$} \put(3.8, 2.1){$2 C_{7}$} \put(7, .5){$1 C_{0}$}
\end{picture}
\end{equation*}%
\begin{equation*}
\underset{\text{Figure 1. Dynkin diagrams of affine }ADE\text{ types}}{%
\setlength{\unitlength}{0.8cm}\begin{picture}(6, 3) \put(-2,
1){$\hat{E}_8:$}\put(0, 1){\circle*{.2}} \put(1.2, 1){\circle*{.2}}
\put(2.4, 1){\circle*{.2}} \put(3.6, 1){\circle*{.2}} \put(4.8,
1){\circle*{.2}} \put(6, 1){\circle*{.2}} \put(7.2,1){\circle*{.2}}
\put(8.4,1){\circle*{.2}} \put(6,2.2){\circle*{.2}} \put(0, 1){\line(1,
0){1.2}} \put(1.2, 1){\line(1, 0){1.2}} \put(2.4, 1){\line(1, 0){1.2}}
\put(3.6, 1){\line(1, 0){1.2}} \put(4.8,1){\line(1,0){1.2}}
\put(6,1){\line(1,0){1.2}} \put(7.2,1){\line(1,0){1.2}} \put(6, 1){\line(0,
1){1.2}} \put(-0.5, .5){$1 C_0$} \put(1.1, .5){$2 C_1$} \put(2.3, .5){$3
C_{2}$} \put(3.5, .5){$4 C_{3}$} \put(4.7, .5){$5 C_{4}$} \put(5.8, .5){$6
C_{5}$} \put(7.2, .5){$4 C_{6}$} \put(6.2, 2.1){$3 C_{8}$} \put(8.4, .5){$2
C_{7}$} \end{picture}\ \ }
\end{equation*}

\begin{remark}
We will also call a nodal or cuspidal rational curve with trivial normal
bundle an $\widehat{A}_{0}$ curve.\qquad
\end{remark}

\begin{remark}
By Kodaira's classification of fibers of relative minimal elliptic surfaces,
every singular fiber is an affine $ADE$ curve unless it is rational with a
cusp, tacnode or triplepoint $($corresponding to type $II$ or $III(\widehat{A%
}_{1})$ or $VI(\widehat{A}_{2})$ in Kodaira's notations$)$, which can also
be regarded as a degenerated affine $ADE$ curve of type $\widehat{A}_{0}$, $%
\widehat{A}_{1}$ or $\widehat{A}_{2}$ respectively. In this paper, we will
not distinguish affine $ADE$ curves from their degenerated forms since they
have the same intersection matrices. We also call the affine $ADE$ curves as
Kodaira curves.
\end{remark}

\begin{definition}
A bundle $E$ is called an $ADE$ $($resp. affine $ADE)$ bundle of type $%
\mathfrak{g}$ $($resp. $\widehat{\mathfrak{g}})$ if $E$ has a fiberwise Lie
algebra structure of the corresponding type.
\end{definition}

In the following two subsections, we will recall an explicit construction of
the Lie algebra $\mathfrak{g}$-bundles, loop Lie algebra $L\mathfrak{g}$%
-bundles and the affine Lie algebra $\widehat{\mathfrak{g}}$-bundles from
(affine) $ADE$ curves in $X$.

\subsection{$ADE$ bundles}

Suppose $C=\cup _{i=1}^{r}C_{i}$ is an $ADE$ curve of type $\mathfrak{g}$ in 
$X$, we will construct the corresponding $ADE$ bundle $\mathcal{E}_{0}^{%
\mathfrak{g}}$ over $X$ as follows \cite{CL}.

Note the rank $r$ of $\mathfrak{g}$ equals the number of $C_{i}$'s, we
denote $\Phi :=\{\alpha =[\sum_{i=1}^{r}a_{i}C_{i}]\in H^{2}(X,\mathbb{Z}%
)|\alpha ^{2}=-2\}$, then $\Phi $ is a simply-laced root system of $%
\mathfrak{g}$ with a base $\Delta :=\{[C_{i}]|i=1,2,\cdots ,r\}$. We have a
decomposition $\Phi =\Phi ^{+}\cup \Phi ^{-}$ into positive and negative
roots. We define a bundle $\mathcal{E}_{0}^{(\mathfrak{g},\Phi )}$ over $X$
as follows:%
\begin{equation*}
\mathcal{E}_{0}^{(\mathfrak{g},\Phi )}:=O^{\oplus r}\oplus \bigoplus_{\alpha
\in \Phi }O(\alpha )\text{.}
\end{equation*}%
Here $O(\alpha )=O(\sum_{i=1}^{r}a_{i}C_{i})$ where $\alpha
=[\sum_{i=1}^{r}a_{i}C_{i}]$. There is an inner product $\langle ,\rangle $
on $\Phi $ defined by $\langle \alpha ,\beta \rangle :=-\alpha \cdot \beta $%
, negative of the intersection form.

For every open chart $U$ of $X$, we take $x_{\alpha }^{U}$ to be a
nonvanishing section of $O_{U}(\alpha )$ and $h_{i}^{U}$ ($1\leq i\leq r$)
nonvanishing sections of $O_{U}^{\oplus r}$. Define a Lie algebra structure $%
[$ $,$ $]_{\Phi }$ on $\mathcal{E}_{0}^{(\mathfrak{g},\Phi )}$ such that $%
\{x_{\alpha }^{U}$'s, $h_{i}^{U}$'s$\}$ is the Chevalley basis \cite{H}, i.e.

(a) $[h_{i}^{U},$ $h_{j}^{U}]_{\Phi }=0$, $1\leq i$, $j\leq r$.

(b) $[h_{i}^{U},$ $x_{\alpha }^{U}]_{\Phi }=\left\langle \alpha ,\text{ }%
C_{i}\right\rangle x_{\alpha }^{U}$, $1\leq i\leq r$, $\alpha \in \Phi $.

(c) $[x_{\alpha }^{U},$ $x_{-\alpha }^{U}]_{\Phi }=h_{\alpha }^{U}$ is a $%
\mathbb{Z}$-linear combination of $h_{i}^{U}$.

(d) If $\alpha $, $\beta $ are independent roots, and $\beta -p\alpha $, $%
\cdots $, $\beta +q\alpha $ is the $\alpha $-string through $\beta $, then $%
[x_{\alpha }^{U},$ $x_{\beta }^{U}]_{\Phi }=0$ if $q=0$, otherwise $%
[x_{\alpha }^{U},$ $x_{\beta }^{U}]_{\Phi }=\pm (p+1)x_{\alpha +\beta }^{U}$.

Since $\mathfrak{g}$ is a simply-laced Lie algebra, all the roots for $%
\mathfrak{g}$ have the same length, we have any $\alpha $-string through $%
\beta $ is of length at most $2$. So (d) can be written as $[x_{\alpha
}^{U}, $ $x_{\beta }^{U}]_{\Phi }=n_{\alpha ,\beta }x_{\alpha +\beta }^{U}$,
where $n_{\alpha ,\beta }=\pm 1$ if $\alpha +\beta \in \Phi $, otherwise $%
n_{\alpha ,\beta }=0$. It is easy to check that these Lie algebra structures
are compatible with different trivializations of $\mathcal{E}_{0}^{(%
\mathfrak{g},\Phi )}$. Hence $\mathcal{E}_{0}^{(\mathfrak{g},\Phi )}$ is a
Lie algebra bundle of type $\mathfrak{g}$ over $X$.

\subsection{Affine $ADE$ bundles}

Suppose $C=\cup _{i=0}^{r}C_{i}$ is an affine $ADE$ curve of type $\widehat{%
\mathfrak{g}}$ in $X$, we will construct the corresponding affine $ADE$
bundle $\mathcal{E}_{0}^{\widehat{\mathfrak{g}}}$ of type $\widehat{%
\mathfrak{g}}$ over $X$ as follows.

First, we choose an extended root of $\widehat{\mathfrak{g}}$, say $C_{0}$,
then $\mathfrak{g}$ is corresponding to the Dynkin diagram consists of those 
$C_{i}$ with $i\neq 0$, i.e. $\Phi :=\{\alpha =[\sum_{i\neq 0}a_{i}C_{i}]\in
H^{2}(X,\mathbb{Z})|\alpha ^{2}=-2\}$ is the root system of $\mathfrak{g}$.
As above, we have a $\mathfrak{g}$-bundle $\mathcal{E}_{0}^{(\mathfrak{g}%
,\Phi )}=O^{\oplus r}\oplus \bigoplus_{\alpha \in \Phi }O(\alpha )$. We
define%
\begin{equation*}
\mathcal{E}_{0}^{(L\mathfrak{g},\Phi )}:=\bigoplus_{n\in \mathbb{Z}}(%
\mathcal{E}_{0}^{(\mathfrak{g},\Phi )}\otimes O(nF))\text{ and }\mathcal{E}%
_{0}^{(\widehat{\mathfrak{g}},\Phi )}:=\bigoplus_{n\in \mathbb{Z}}(\mathcal{E%
}_{0}^{(\mathfrak{g},\Phi )}\otimes O(nF))\oplus O\text{.}
\end{equation*}

We know $\Phi _{\widehat{\mathfrak{g}}}:=\{\alpha +nF|\alpha \in \Phi ,n\in 
\mathbb{Z}\}\cup \{nF|n\in \mathbb{Z},n\neq 0\}$ is an affine root system
and it decomposes into union of positive and negative roots, i.e. $\Phi _{%
\widehat{\mathfrak{g}}}=\Phi _{\widehat{\mathfrak{g}}}^{+}\cup \Phi _{%
\widehat{\mathfrak{g}}}^{-}$, where $\Phi _{\widehat{\mathfrak{g}}%
}^{+}=\{\sum a_{i}C_{i}\in \Phi _{\widehat{\mathfrak{g}}}|a_{i}\geq 0$ for
all $i\}=\{\alpha +nF|\alpha \in \Phi ^{+},n\in \mathbb{Z}_{\geq 0}\}\cup
\{\alpha +nF|\alpha \in \Phi ^{-},n\in \mathbb{Z}_{\geq 1}\}\cup \{nF|n\in 
\mathbb{Z}_{\geq 1}\}$ and $\Phi _{\widehat{\mathfrak{g}}}^{-}=-\Phi _{%
\widehat{\mathfrak{g}}}^{+}$.

To describe the Lie algebra structures, we proceed as before, for every open
chart $U$ of $X$, we take a local basis $e_{i}^{U}$ of $\mathcal{E}_{0}^{(%
\mathfrak{g},\Phi )}|_{U}$ ($e_{i}^{U}$ is just $h_{j}^{U}$ or $x_{\alpha
}^{U}$ as above), $e_{nF}^{U}$ of $O(nF)|_{U}$, $e_{c}^{U}$ of $O|_{U}$,
compatible with the tensor product, for example, $e_{nF}^{U}\otimes
e_{mF}^{U}=e_{(n+m)F}^{U}$. Then define%
\begin{equation}
\lbrack e_{i}^{U}e_{nF}^{U},e_{j}^{U}e_{mF}^{U}]_{L\mathfrak{g},\Phi
}:=[e_{i}^{U},e_{j}^{U}]_{\Phi }e_{(n+m)F}^{U}\text{,}  \label{1}
\end{equation}%
\begin{equation}
\lbrack e_{i}^{U}e_{nF}^{U}+\lambda e_{c}^{U},e_{j}^{U}e_{mF}^{U}+\mu
e_{c}^{U}]_{\widehat{\mathfrak{g}},\Phi }:=[e_{i}^{U},e_{j}^{U}]_{\Phi
}e_{(n+m)F}^{U}+n\delta _{n+m,0}k(e_{i}^{U},e_{j}^{U})e_{c}^{U}\text{.}
\label{2}
\end{equation}%
Here $[$ $,$ $]_{\Phi }$ is the Lie bracket on $\mathcal{E}_{0}^{(\mathfrak{g%
},\Phi )}$ and $k(x,y)=Tr(adx$ $ady)$ is the Killing form on $\mathfrak{g}$.

\begin{lemma}
$(1)$ $($resp. $(2))$ defines a fiberwise loop $($resp. affine$)$ Lie
algebra structure which is compatible with any trivialization of $\mathcal{E}%
_{0}^{(L\mathfrak{g},\Phi )}$ $($resp. $\mathcal{E}_{0}^{(\widehat{\mathfrak{%
g}},\Phi )})$.

\begin{proof}
See Proposition 23 of \cite{LXZ}.
\end{proof}
\end{lemma}

From the above lemma, we have the following result.

\begin{proposition}
If $C$ is an affine $ADE$ curve of type $\widehat{\mathfrak{g}}$ in $X$,
then $\mathcal{E}_{0}^{(L\mathfrak{g},\Phi )}$ $($resp. $\mathcal{E}_{0}^{(%
\widehat{\mathfrak{g}},\Phi )})$ is a loop $($resp. affine$)$ Lie algebra
bundle of type $L\mathfrak{g}$ $($resp. $\widehat{\mathfrak{g}})$ over $X$.
\end{proposition}

Note any $C_{i}$ with $n_{i}=1$ can be chosen as the extended root
(Appendix).

\begin{proposition}
\label{affine root}The loop Lie algebra bundle $(\mathcal{E}_{0}^{(L%
\mathfrak{g},\Phi )},[$ $,$ $]_{L\mathfrak{g},\Phi })$ does not depend on
the choice of the extended root.

\begin{proof}
Suppose $C_{k}$ $(k\neq 0)$ is another root with $n_{k}=1$, we denote $\Psi
=\{\beta =[\sum_{i\neq k}b_{i}C_{i}]\in H^{2}(X,\mathbb{Z})|\beta ^{2}=-2\}$%
, then $\Psi $ is a root system of $\mathfrak{g}$. As before, we construct
the Lie algebra bundle $\mathcal{E}_{0}^{(\mathfrak{g,}\Psi )}$ and $%
\mathcal{E}_{0}^{(L\mathfrak{g},\Psi )}$ from $\Psi $.

We denote $\alpha _{0}:=\sum_{i\neq 0}n_{i}C_{i}=F-C_{0}$, the longest root
in $\Phi $. For any $\alpha =\sum_{i\neq 0}a_{i}(\alpha )C_{i}\in \Phi $, $%
a_{k}(\alpha )$ can only be $0$, $\pm 1$. Hence there is a bijection between 
$\Phi $ and $\Psi $ given by $\alpha \mapsto \beta =\alpha -a_{k}(\alpha )F$%
. Then from the definitions of $\mathcal{E}_{0}^{(L\mathfrak{g},\Phi )}$ and 
$\mathcal{E}_{0}^{(L\mathfrak{g},\Psi )}$, we know they are the same as
holomorphic vector bundles.

We compare the Lie brackets on them. We choose a local basis of $\mathcal{E}%
_{0}^{(L\mathfrak{g},\Psi )}$ compatible with those of $\mathcal{E}_{0}^{(L%
\mathfrak{g},\Phi )}$ and define $[,]_{L\mathfrak{g},\Psi }$ similarly as $%
[,]_{L\mathfrak{g},\Phi }$, i.e.

$(i)$ when $\beta =\alpha \in \Phi \cap \Psi $, we take $x_{\beta
}=x_{\alpha }$;

$(ii)$ when $\beta =\alpha +F\in \Psi ^{+}\backslash \Phi $, we take $%
x_{\beta }=x_{\alpha }e_{F}$;

$(iii)$ when $\beta =\alpha -F\in \Psi ^{-}\backslash \Phi $, we take $%
x_{\beta }=x_{\alpha }e_{-F}$;

$(iv)$ take $h_{i}$ $(i\neq 0,k)$ as before, take $h_{0}=-h_{\alpha _{0}}$
as we want $[x_{C_{0}},x_{-C_{0}}]_{L\mathfrak{g},\Psi }=[x_{-\alpha
_{0}+F},x_{\alpha _{0}-F}]_{L\mathfrak{g},\Phi }$.

It is obvious $[$ $,$ $]_{L\mathfrak{g},\Psi }=[$ $,$ $]_{L\mathfrak{g},\Phi
}$ on $\mathcal{E}_{0}^{(L\mathfrak{g},\Psi )}\cong \mathcal{E}_{0}^{(L%
\mathfrak{g},\Phi )}$.
\end{proof}
\end{proposition}

For the affine case, we recall that the Killing form of $\mathfrak{g}$ is
the symmetric bilinear map $k:\mathfrak{g\times g\rightarrow }\mathbb{C}$
defined by $k(x,y)=Tr(adx$ $ady)$. It is $ad$-invariant, that is for $%
x,y,z\in \mathfrak{g}$, $k([x,y],z)=k(x,[y,z])$.

\begin{lemma}
For any simple simply-laced Lie algebra $\mathfrak{g}$ with a Chavelly basis 
$\{x_{\alpha },\alpha \in \Phi ;h_{i},1\leq i\leq r\}$ and $m^{\ast }(%
\mathfrak{g})$ the dual Coxeter number of $\mathfrak{g}$, we have

$(i)$ $k(h_{i},x_{\alpha })=0$ for any $i$ and $\alpha $;

$(ii)$ $k(x_{\alpha },x_{\beta })=0$ for any $\alpha +\beta \neq 0$;

$(iii)$ $k(h_{i},h_{j})=2m^{\ast }(\mathfrak{g})\langle C_{i},C_{j}\rangle $;

$(iv)$ $k(x_{\alpha },x_{-\alpha })=2m^{\ast }(\mathfrak{g})$ for any $%
\alpha $.

\begin{proof}
Directly from the Killing form $k$ being $ad$-invariant or see \cite{Ma}.
\end{proof}
\end{lemma}

\begin{proposition}
The affine Lie algebra bundle $(\mathcal{E}_{0}^{(\widehat{\mathfrak{g}}%
,\Phi )},[$ $,$ $]_{\widehat{\mathfrak{g}},\Phi })$ does not depend on the
choice of the extended root.

\begin{proof}
Follow the notations in Proposition \ref{affine root}, but we will take $%
h_{0}=-h_{\alpha _{0}}+2m^{\ast }(\mathfrak{g})e_{c}$. We will check that $[$
$,$ $]_{\widehat{\mathfrak{g}},\Psi }=[$ $,$ $]_{\widehat{\mathfrak{g}},\Phi
}$ on $\mathcal{E}_{0}^{(\widehat{\mathfrak{g}},\Psi )}=\mathcal{E}_{0}^{(%
\widehat{\mathfrak{g}},\Phi )}$:

$(a)$ when $\beta _{1}=\alpha _{1}+F,\beta _{2}=\alpha _{2}+F\in \Psi
^{+}\backslash \Phi $, $\alpha _{1},\alpha _{2}\in \Phi ^{-}\backslash \Psi $
we have%
\begin{equation*}
\lbrack h_{\beta _{1}}e_{nF},h_{\beta _{2}}e_{mF}]_{\widehat{\mathfrak{g}}%
,\Psi }=n\delta _{n+m,0}k(h_{\beta _{1}},h_{\beta _{2}})e_{c}\text{,}
\end{equation*}%
which is the same with%
\begin{equation*}
\lbrack h_{-\alpha _{1}}e_{nF},h_{-\alpha _{2}}e_{mF}]_{\widehat{\mathfrak{g}%
},\Phi }=n\delta _{n+m,0}k(h_{\alpha _{1}},h_{\alpha _{2}})e_{c}\text{,}
\end{equation*}%
since $k(h_{\beta _{1}},h_{\beta _{2}})=2m^{\ast }(\mathfrak{g})\langle
\beta _{1},\beta _{2}\rangle =2m^{\ast }(\mathfrak{g})\langle F-\alpha
_{1},F-\alpha _{2}\rangle =k(h_{\alpha _{1}},h_{\alpha _{2}})$.

$(b)$ For $[h_{i}e_{nF},x_{\alpha }e_{mF}]_{\widehat{\mathfrak{g}},\Phi }$,
automatically from $k(h_{i},x_{\alpha })=0$ and loop case.

$(c)$ When $\beta =\alpha +F\in \Psi ^{+}\backslash \Phi $, $\alpha \in \Phi
^{-}\backslash \Psi $,%
\begin{equation*}
\lbrack x_{\beta }e_{nF},x_{-\beta }e_{mF}]_{\widehat{\mathfrak{g}},\Psi
}=h_{\beta }e_{(n+m)F}+n\delta _{n+m,0}k(x_{\beta },x_{-\beta })e_{c}\text{,}
\end{equation*}%
which is the same with%
\begin{equation*}
\lbrack x_{-\alpha }e_{(n+1)F},x_{\alpha }e_{(m-1)F}]_{\widehat{\mathfrak{g}}%
,\Phi }=-h_{\alpha }e_{(n+m)F}+(n+1)\delta _{n+m,0}k(x_{\alpha },x_{-\alpha
})e_{c}\text{,}
\end{equation*}%
by considering $m+n=0$ and $m+n\neq 0$ separately.

$(d)$ For $[x_{\alpha _{1}}e_{nF},x_{\alpha _{2}}e_{mF}]_{\widehat{\mathfrak{%
g}},\Phi }$ with $\alpha _{1}+\alpha _{2}\neq 0$, automatically from $%
k(x_{\alpha _{1}},x_{\alpha _{2}})=0$ and loop case.
\end{proof}
\end{proposition}

For simplicity, we will omit $\Phi $ in $(\mathfrak{g,}\Phi )$, $(L\mathfrak{%
g},\Phi )$ and $(\widehat{\mathfrak{g}},\Phi )$ when there is no confusion.

\section{Trivialization of $\mathcal{E}_{0}^{\widehat{\mathfrak{g}}}$ over $%
C_{i}$'s after deformations}

If $C=\cup C_{i}$ is an affine $ADE$ curve in $X$, then the corresponding $%
F=\sum n_{i}C_{i}$ satisfies $F\cdot F=0$, i.e. $O_{F}(F)$ is a
topologically trivial bundle. If $O_{F}(F)$ is trivial holomorphically and $%
q(X)=0$, then from the long exact sequence of cohomologies induced by $%
0\rightarrow O_{X}\rightarrow O_{X}(F)\rightarrow O_{F}(F)\rightarrow 0$, we
know $H^{0}(X,O_{X}(F))\cong \mathbb{C}^{2}$. Hence $F$ is a fiber of an
elliptic fibration on $X$.

Suppose $X$ is an elliptic surface, i.e. there is a smooth curve $B$ and a
surjective morphism $\pi :X\rightarrow B$ whose generic fiber $F_{b}$ ($b\in
B$) is an elliptic curve. Assume $\pi $ is singular at $b_{0}\in B$ and $%
F_{b_{0}}=\sum n_{i}C_{i}$ is a singular fiber of type $\widehat{\mathfrak{g}%
}$. Hence, we have a $\widehat{\mathfrak{g}}$-bundle $\mathcal{E}_{0}^{%
\widehat{\mathfrak{g}}}$ over $X$. The restriction of $\mathcal{E}_{0}^{%
\widehat{\mathfrak{g}}}$ to any fiber $F_{b}$, other than $F_{b_{0}}$, is
trivial because $F_{b}\cap C_{i}=\varnothing $ for any $i$. However, $%
\mathcal{E}_{0}^{\widehat{\mathfrak{g}}}|_{F_{b_{0}}}$ is not trivial, for
instance $O(-C_{i})|_{C_{i}}\cong O_{\mathbb{P}^{1}}(2)$. Nevertheless, we
will show that after deformations of holomorphic structures, $\mathcal{E}%
_{0}^{\widehat{\mathfrak{g}}}$ will become trivial on every irreducible
component of $F_{b_{0}}$.

\subsection{Review of $ADE$ cases}

In our earlier paper \cite{CL}, we showed how to take successive extensions
to make the $\mathfrak{g}$-bundle $\mathcal{E}_{0}^{\mathfrak{g}}$ trivial
on every component $C_{i}$ of the $ADE$ curve $C=\cup _{i=1}^{r}C_{i}$.

\begin{definition}
Given any $\varphi =(\varphi _{\alpha })_{\alpha \in \Phi ^{+}}\in \Omega
^{0,1}(X,\bigoplus_{\alpha \in \Phi ^{+}}O(\alpha ))$, we define $\overline{%
\partial }_{\varphi }:\Omega ^{0,0}(X,\mathcal{E}_{0}^{\mathfrak{g}%
})\longrightarrow \Omega ^{0,1}(X,\mathcal{E}_{0}^{\mathfrak{g}})$ by 
\begin{equation*}
\overline{\partial }_{\varphi }:=\overline{\partial }_{0}+ad(\varphi ):=%
\overline{\partial }_{0}+\sum_{\alpha \in \Phi ^{+}}ad(\varphi _{\alpha })%
\text{,}
\end{equation*}
\end{definition}

More explicitly, if we write $\varphi _{\alpha }=c_{\alpha }^{U}x_{\alpha
}^{U}$ locally for some one form $c_{\alpha }^{U}$, then $ad(\varphi
_{\alpha })=c_{\alpha }^{U}ad(x_{\alpha }^{U})$. It is easy to check that $%
\overline{\partial }_{\mathcal{\varphi }}$ is well-defined and compatible
with the Lie algebra structure, i.e. $\overline{\partial }_{\mathcal{\varphi 
}}[,]_{\Phi }=0$. For $\overline{\partial }_{\mathcal{\varphi }}$ to define
a holomorphic structure, we need 
\begin{equation*}
0=\overline{\partial }_{\mathcal{\varphi }}^{2}=\sum_{\alpha \in \Phi ^{+}}(%
\overline{\partial }_{0}c_{\alpha }^{U}+\sum_{\beta +\gamma =\alpha
}(n_{\beta ,\gamma }c_{\beta }^{U}c_{\gamma }^{U}))ad(x_{\alpha }^{U})\text{,%
}
\end{equation*}%
that is $\overline{\partial }_{0}\varphi _{\alpha }+\sum_{\beta +\gamma
=\alpha }(n_{\beta ,\gamma }\varphi _{\beta }\varphi _{\gamma })=0$ for any $%
\alpha \in \Phi ^{+}$. Explicitly:%
\begin{equation*}
\left\{ 
\begin{tabular}{cc}
$\overline{\partial }_{0}\varphi _{C_{i}}=0$ & $i=1,2\cdots ,r$ \\ 
$\overline{\partial }_{0}\varphi _{C_{i}+C_{j}}=n_{C_{i},C_{j}}\varphi
_{C_{i}}\varphi _{C_{j}}$ & $\text{ if }C_{i}+C_{j}\in \Phi ^{+}$ \\ 
$\vdots $ & 
\end{tabular}%
\ \ \right.
\end{equation*}%
Recall $\{C_{i}\}_{i=1}^{r}\subset \Phi ^{+}$ is a base.

\begin{proposition}
\label{holomorphic}Given any $(\varphi _{C_{i}})_{i=1}^{r}\in \Omega
^{0,1}(X,\bigoplus_{i=1}^{r}O(C_{i}))$ with $\overline{\partial }\varphi
_{C_{i}}=0$ for any $i$, it can be extended to $\varphi =(\varphi _{\alpha
})_{\alpha \in \Phi ^{+}}\in \Omega ^{0,1}(X,\bigoplus_{\alpha \in \Phi
^{+}}O(\alpha ))$ satisfying $\overline{\partial }_{\mathcal{\varphi }%
}^{2}=0 $. Namely we have a holomorphic $\mathfrak{g}$-bundle $\mathcal{E}%
_{\varphi }^{\mathfrak{g}}$ over $X$.
\end{proposition}

The proof of this proposition uses the following lemma.

\begin{lemma}
\label{cohomology}If $\mathit{p}_{g}(X)=0$, then

$(i)$ for any $\alpha \in \Phi ^{+}$, $H^{2}(X,O(\alpha ))=0$.

$(ii)$ the restriction homomorphism $H^{1}(X,$ $O_{X}(C_{i}))\rightarrow
H^{1}(X,$ $O_{C_{i}}(C_{i}))$ is surjective.
\end{lemma}

\begin{theorem}
\label{thm1}For any given $i$, the holomorphic $\mathfrak{g}$-bundle $%
\mathcal{E}_{\varphi }^{\mathfrak{g}}$ over $X$ is trivial on $C_{i}$ if and
only if $[\varphi _{C_{i}}|_{C_{i}}]\neq 0$.
\end{theorem}

Note that part (ii) of Lemma \ref{cohomology} says that such $\varphi
_{C_{i}}$'s can always be found.

\subsection{Trivializations in loop $ADE$ cases}

\begin{definition}
Given any $\varphi =(\varphi _{\alpha })_{\alpha \in \Phi _{\widehat{%
\mathfrak{g}}}^{+}}\in \Omega ^{0,1}(X,\bigoplus_{\alpha \in \Phi _{\widehat{%
\mathfrak{g}}}^{+}}O(\alpha ))$, we define $\overline{\partial }_{(\varphi
,\Phi )}:\Omega ^{0,0}(X,\mathcal{E}_{0}^{L\mathfrak{g}})\longrightarrow
\Omega ^{0,1}(X,\mathcal{E}_{0}^{L\mathfrak{g}})$ by $\overline{\partial }%
_{(\varphi ,\Phi )}:=\overline{\partial }_{0}+ad(\varphi )$.
\end{definition}

More explicitly, similarly as explained in $\S 3.1$, we have 
\begin{eqnarray*}
\overline{\partial }_{(\varphi ,\Phi )} &:&=\overline{\partial }%
_{0}+\sum_{n\in \mathbb{Z}_{\geq 0}}\sum_{\alpha \in \Phi ^{+}}(c_{\alpha
+nF}e_{nF}ad(x_{\alpha })+c_{-\alpha +(n+1)F}e_{(n+1)F}ad(x_{-\alpha })) \\
&&+\sum_{n\in \mathbb{Z}_{\geq
0}}\sum_{i=1}^{r}c_{(n+1)F}^{i}e_{(n+1)F}ad(h_{i})\text{,}
\end{eqnarray*}

\begin{proposition}
\label{Lie bracket2}$\overline{\partial }_{(\varphi ,\Phi )}$ is compatible
with the Lie algebra structure on $\mathcal{E}_{0}^{L\mathfrak{g}}$.

\begin{proof}
$\overline{\partial }_{(\varphi ,\Phi )}[$ $,$ $]_{L\mathfrak{g},\Phi }=0$
follows directly from the Jacobi identity.
\end{proof}
\end{proposition}

For $\overline{\partial }_{(\varphi ,\Phi )}$ to define a holomorphic
structure, we need $\overline{\partial }_{(\varphi ,\Phi )}^{2}=0$, which is
equivalent to the following equations:%
\begin{equation*}
\left\{ 
\begin{array}{l}
\overline{\partial }_{0}\varphi _{nF}^{i}=\sum_{p+q=n}\sum_{\alpha \in \Phi
^{+}}\pm a_{i}(h_{\alpha })\varphi _{\alpha +pF}\varphi _{-\alpha +qF}\text{,%
} \\ 
\\ 
\overline{\partial }_{0}\varphi _{\alpha +nF}=\sum_{p+q=n}(\sum_{\alpha
_{1}+\alpha _{2}=\alpha }\pm \varphi _{\alpha _{1}+pF}\varphi _{\alpha
_{2}+qF}+\sum_{i=1}^{r}\left\langle \alpha ,C_{i}\right\rangle \varphi
_{\alpha +pF}\varphi _{qF}^{i})\text{,} \\ 
\\ 
\overline{\partial }_{0}\varphi _{-\alpha +nF}=\sum_{p+q=n}(\sum_{\alpha
_{2}-\alpha _{1}=\alpha }\pm \varphi _{\alpha _{1}+pF}\varphi _{-\alpha
_{2}+qF}+\sum_{i=1}^{r}\left\langle -\alpha ,C_{i}\right\rangle \varphi
_{-\alpha +pF}\varphi _{qF}^{i})\text{,}%
\end{array}%
\right.
\end{equation*}%
where $a_{i}(h_{\alpha })$ is the coefficient of $h_{i}$ in $h_{\alpha }$.

\begin{proposition}
\label{holomorphic2}Given any $(\varphi _{C_{i}})_{i=0}^{r}\in \Omega
^{0,1}(X,\bigoplus_{i=0}^{r}O(C_{i}))$ with $\overline{\partial }\varphi
_{C_{i}}=0$ for every $i$, it can be extended to $\varphi =(\varphi _{\alpha
})_{\alpha \in \Phi _{\widehat{\mathfrak{g}}}^{+}}\in \Omega
^{0,1}(X,\bigoplus_{\alpha \in \Phi _{\widehat{\mathfrak{g}}}^{+}}O(\alpha
)) $ satisfying $\overline{\partial }_{\mathcal{\varphi }}^{2}=0$. Namely we
have a holomorphic $L\mathfrak{g}$-bundle $\mathcal{E}_{\varphi }^{L%
\mathfrak{g}}$ over $X$.
\end{proposition}

In order to prove this proposition, we need the following lemma.

\begin{lemma}
\label{cohomology2}If $\mathit{p}_{g}(X)=0$, then for any $\alpha \in \Phi
^{+}$, $n\in \mathbb{Z}_{\geq 0}$, $H^{2}(X,O(nF))$, $H^{2}(X,O(\alpha +nF))$
and $H^{2}(X,O(-\alpha +(n+1)F))$ are zeros.

\begin{proof}
Since $F$ is an effective divisor and $H^{0}(X,K_{X})=0$, we have for any $%
n\geq 0$, $H^{0}(X,K_{X}(-nF))=0$. This is equivalent to $H^{2}(X,O(nF))=0$
by Serre duality. Similarly, $H^{2}(X,O(\alpha +nF))=0$ follows from $%
H^{0}(X,K_{X}(-\alpha ))\cong H^{2}(X,O(\alpha ))=0$ $($Lemma \ref%
{cohomology}$)$. The proof of $H^{2}(X,O(-\alpha +(n+1)F))=0$ uses the fact
that $F-\alpha $ is an effective divisor for any $\alpha \in \Phi ^{+}$.
\end{proof}
\end{lemma}

\begin{proof}
$\mathit{(}$\textit{of Proposition \ref{holomorphic2}}$\mathit{)}$\textit{:
the equation }$\overline{\partial }_{(\varphi ,\Phi )}^{2}=0$\textit{\ can
be rewritten as follows:}%
\begin{equation*}
\left\{ 
\begin{tabular}{l}
$\overline{\partial }_{0}\varphi _{C_{i}}=0\text{ for }i=1,2\cdots ,r\text{,}
$ \\ 
$\overline{\partial }_{0}\varphi _{\alpha }=\sum_{\alpha _{1}+\alpha
_{2}=\alpha }(\pm \varphi _{\alpha _{1}}\varphi _{\alpha _{2}})\text{,}$ \\ 
$\overline{\partial }_{0}\varphi _{-\alpha _{0}+F}=\overline{\partial }%
_{0}\varphi _{C_{0}}=0\text{,}$ \\ 
$\overline{\partial }_{0}\varphi _{-\alpha +F}=\sum_{\alpha _{2}-\alpha
_{1}=\alpha }(\pm \varphi _{\alpha _{1}}\varphi _{-\alpha _{2}+F})\text{,}$
\\ 
$\overline{\partial }_{0}\varphi _{F}^{i}=\sum_{\alpha \in \Phi ^{+}}(\pm
a_{i}(h_{\alpha })\varphi _{\alpha }\varphi _{-\alpha +F})\text{,}$ \\ 
$\vdots $%
\end{tabular}%
\right.
\end{equation*}
\textit{where }$\alpha _{0}=F-C_{0}$\textit{\ is the longest root in }$\Phi $%
\textit{.}

\textit{Firstly, we can solve for all the }$\varphi _{\alpha }$'s\textit{, }$%
\alpha \in \Phi ^{+}$\textit{\ from }$H^{2}(X,O(\alpha ))=0$ $($\textit{%
Proposition \ref{holomorphic}}$)$\textit{. Secondly, we get all the }$%
\varphi _{-\alpha +F}$'s\textit{, }$\alpha \in \Phi ^{+}$\textit{\ from }$%
H^{2}(X,O(-\alpha +F))=0$\textit{. Thirdly, since we have all the }$\varphi
_{\alpha }$\textit{'s and }$\varphi _{-\alpha +F}$\textit{'s, we can solve
for all the }$\varphi _{F}^{i}$'s\textit{\ for }$1\leq i\leq r$\textit{\
from }$H^{2}(X,O(F))=0$\textit{. Do this process for }$\varphi _{\alpha +nF}$%
\textit{, }$\varphi _{-\alpha +(n+1)F}$\textit{\ and }$\varphi _{(n+1)F}^{i}$%
\textit{\ inductively on }$n$\textit{.}
\end{proof}

\bigskip

By Lemma \ref{cohomology}, there always exists $\varphi _{C_{i}}\in \Omega
^{0,1}(X,$ $O(C_{i}))$ such that $0\neq \lbrack \varphi
_{C_{i}}|_{C_{i}}]\in H^{1}(X,$ $O_{C_{i}}(C_{i}))\cong \mathbb{C}$ for each 
$i=0,1,\cdots r$.

\begin{theorem}
\label{loop thm}For any given $i$, the holomorphic $L\mathfrak{g}$-bundle $%
\mathcal{E}_{\varphi }^{L\mathfrak{g}}$ over $X$ is trivial on $C_{i}$ if
and only if $[\varphi _{C_{i}}|_{C_{i}}]\neq 0$.

\begin{proof}
The proof will be given in $\S 3.4$ and $\S 3.5$. In $\S 3.4$, we deal with
all the loop $ADE$ cases except loop $E_{8}$ case which will be analyzed in $%
\S 3.5$.
\end{proof}
\end{theorem}

\subsection{Trivializations in loop $ADE$ cases}

Follow the notations in $\S 3.2$, we define $\overline{\partial }_{(\varphi
,\Phi )}:=\overline{\partial }_{0}+ad(\varphi )$ on $\mathcal{E}_{0}^{%
\widehat{\mathfrak{g}}}$, note the adjoint action here is defined using the
affine Lie bracket.

\begin{proposition}
\label{Lie bracket3}$\overline{\partial }_{(\varphi ,\Phi )}$ is compatible
with the Lie algebra structure on $\mathcal{E}_{0}^{\widehat{\mathfrak{g}}}$.

\begin{proof}
$\overline{\partial }_{(\varphi ,\Phi )}[$ $,$ $]_{\widehat{\mathfrak{g}}%
,\Phi }=0$ follows directly from the Jacobi identity and the Killing from
being invariant under the adjoint action.
\end{proof}
\end{proposition}

It is easy to see that $\overline{\partial }_{(\varphi ,\Phi )}^{2}=0$ in
the affine case is equivalent to $\overline{\partial }_{(\varphi ,\Phi
)}^{2}=0$ in the loop case. Hence we have a new holomorphic structure $%
\overline{\partial }_{(\varphi ,\Phi )}$ on $\mathcal{E}_{0}^{\widehat{%
\mathfrak{g}}}$.

\begin{theorem}
\label{thm2}For any given $i$, the holomorphic $\widehat{\mathfrak{g}}$%
-bundle $\mathcal{E}_{\varphi }^{\widehat{\mathfrak{g}}}$ over $X$ is
trivial on $C_{i}$ if and only if $[\varphi _{C_{i}}|_{C_{i}}]\neq 0$.

\begin{proof}
\textit{This follows from Theorem \ref{loop thm},\ }$0\rightarrow
O\rightarrow E_{\varphi }^{\widehat{\mathfrak{g}}}\rightarrow E_{\varphi }^{L%
\mathfrak{g}}\rightarrow 0$ and $Ext_{\mathbb{P}^{1}}^{1}(O,O)=0$\textit{.}
\end{proof}
\end{theorem}

\subsection{Proof (except the loop $E_{8}$ case)}

In this subsection, we use the symmetry of the affine $ADE$ Dynkin diagram
(except $\widehat{E}_{8}$) to show that $\mathcal{E}_{\varphi }^{L\mathfrak{g%
}}$ is trivial on $C_{i}$ if an only if $[\varphi _{C_{i}}|_{C_{i}}]\neq 0$.

Recall, topologically, $\mathcal{E}_{\varphi }^{L\mathfrak{g}}$ is $\mathcal{%
E}_{0}^{L\mathfrak{g}}=\bigoplus_{n\in \mathbb{Z}}(\mathcal{E}_{0}^{(%
\mathfrak{g},\Phi )}\otimes O(nF))$, but with a holomorphic structure $%
\overline{\partial }_{(\varphi ,\Phi )}$ of the following upper triangular
block shape:%
\begin{equation*}
\overline{\partial }_{\varphi }=\left( 
\begin{tabular}{c|c|c|c|c}
$\ddots $ & $\ddots $ & $\ddots $ & $\ddots $ & $\ddots $ \\ \hline
$\ddots $ & $\overline{\partial }_{\mathcal{E}_{\varphi }^{(\mathfrak{g}%
,\Phi )}\otimes O((n+1)F)}$ & $\ast $ & $\ast $ & $\ddots $ \\ \hline
$\ddots $ & $0$ & $\overline{\partial }_{\mathcal{E}_{\varphi }^{(\mathfrak{g%
},\Phi )}\otimes O(nF)}$ & $\ast $ & $\ddots $ \\ \hline
$\ddots $ & $0$ & $0$ & $\overline{\partial }_{\mathcal{E}_{\varphi }^{(%
\mathfrak{g},\Phi )}\otimes O((n-1)F)}$ & $\ddots $ \\ \hline
$\ddots $ & $\ddots $ & $\ddots $ & $\ddots $ & $\ddots $%
\end{tabular}%
\ \right) \text{.}
\end{equation*}%
i.e. $\mathcal{E}_{\varphi }^{L\mathfrak{g}}$ is constructed from successive
extensions of $\mathcal{E}_{\varphi }^{(\mathfrak{g},\Phi )}\otimes O(nF)$'s.

Note $\overline{\partial }_{(\varphi ,\Phi )}|_{\mathcal{E}_{\varphi }^{(%
\mathfrak{g},\Phi )}}=\overline{\partial }_{0}+\sum_{\alpha \in \Phi
^{+}}ad(\varphi _{\alpha })$. By Theorem \ref{thm1}, for every $i\neq 0$, $%
\mathcal{E}_{\varphi }^{(\mathfrak{g},\Phi )}$ is trivial on $C_{i}$ if and
only if $[\varphi _{C_{i}}|_{C_{i}}]\neq 0$. We also know $O(F)|_{C_{i}}$ is
trivial for every $i$ because $F\cdot C_{i}=0$. Thus, when $i\neq 0$, $%
\mathcal{E}_{\varphi }^{L\mathfrak{g}}|_{C_{i}}$ is constructed from
successive extensions of trivial vector bundles over $C_{i}\cong \mathbb{P}%
^{1}$. This implies that $\mathcal{E}_{\varphi }^{L\mathfrak{g}}|_{C_{i}}$
is trivial if and only if $[\varphi _{C_{i}}|_{C_{i}}]\neq 0$ as $Ext_{%
\mathbb{P}^{1}}^{1}(O,O)=0$.

Now we consider $i=0$. Since $\widehat{\mathfrak{g}}\neq \widehat{E}_{8}$,
the affine Dynkin diagram always admits a diagram automorphism, that means
we can write $\mathcal{E}_{0}^{L\mathfrak{g}}$ as $\bigoplus_{n\in \mathbb{Z}%
}(\mathcal{E}_{0}^{(\mathfrak{g},\Psi )}\otimes O(nF))$ (see Proposition \ref%
{affine root}). Suppose the extended root corresponding to $\Psi $ is $C_{k}$%
, and the longest root in $\Psi $ is $\beta _{0}$.

We will rewrite the holomorphic structure $\overline{\partial }_{(\varphi
,\Phi )}$ in terms of the $\Psi $ root system. Note $\overline{\partial }%
_{(\varphi ,\Phi )}$ is determined by the loop Lie algebra structure which
is independent of the choice of the extended root. We choose a local base of 
$\mathcal{E}_{0}^{(\mathfrak{g},\Psi )}$ as in Proposition \ref{affine root}
and define $\overline{\partial }_{(\psi ,\Psi )}$ to be the same with $%
\overline{\partial }_{(\varphi ,\Phi )}$, then obviously $\psi _{D}=\varphi
_{D}$ when $D\neq nF$.

Because $(\mathcal{E}_{\varphi }^{(L\mathfrak{g},\Phi )},\overline{\partial }%
_{(\varphi ,\Phi )})=(\mathcal{E}_{\psi }^{(L\mathfrak{g},\Psi )},\overline{%
\partial }_{(\psi ,\Psi )})$ as a holomorphic vector bundle, similar to the
arguments in $(\mathcal{E}_{\varphi }^{(L\mathfrak{g},\Phi )},\overline{%
\partial }_{(\varphi ,\Phi )})$ case, we have when $i\neq k$, $\mathcal{E}%
_{\varphi }^{L\mathfrak{g}}$ is trivial on $C_{i}$ if and only if $[\psi
_{C_{i}}|_{C_{i}}]\neq 0$. Note $\psi _{C_{0}}=\varphi _{-\alpha
_{0}+F}=\varphi _{C_{0}}$. So we have Theorem \ref{loop thm} when $\mathfrak{%
g}\neq E_{8}$.

\subsection{Proof for the loop $E_{8}$ case}

Similar to the above subsection, we have when $i=1,2,\cdots 8$, $\mathcal{E}%
_{\varphi }^{LE_{8}}$ is trivial on $C_{i}$ if and only if $[\varphi
_{C_{i}}|_{C_{i}}]\neq 0$. The question is what about $C_{0}$?

We recall $\mathcal{E}_{0}^{E_{8}}:=O^{\oplus 8}\oplus \bigoplus_{\alpha \in
\Phi }O(\alpha )$. For any $\alpha \in \Phi $, we write $a_{1}(\alpha )$ as
the coefficient of $C_{1}$ in $\alpha $, then $O(\alpha )|_{C_{0}}\cong O_{%
\mathbb{P}^{1}}(a_{1}(\alpha ))$. Among $\Phi ^{+}$, there are $63$ roots
with $a_{1}(\alpha )=0$, corresponding to the positive roots of the Lie
sub-algebra $E_{7}$; $56$ roots with $a_{1}(\alpha )=1$, corresponding to
weights of the standard representation of $E_{7}$; $1$ root with $%
a_{1}(\alpha )=2$, which is just the longest root $\alpha _{0}=F-C_{0}$. We
denote $\mathcal{E}_{0}^{E_{7}}\triangleq O^{\oplus 7}\oplus
\bigoplus_{\alpha \in \Phi ,a_{1}(\alpha )=0}O(\alpha )$, $%
V_{0}^{+}\triangleq \bigoplus_{\alpha \in \Phi ,a_{1}(\alpha )=1}O(\alpha )$
and $V_{0}^{-}\triangleq \bigoplus_{\alpha \in \Phi ,a_{1}(\alpha
)=-1}O(\alpha )$, then%
\begin{equation*}
\mathcal{E}_{0}^{E_{8}}=\mathcal{E}_{0}^{E_{7}}\oplus O\oplus
V_{0}^{+}\oplus V_{0}^{-}\oplus O(\alpha _{0})\oplus O(-\alpha _{0})\text{.}
\end{equation*}

When $O(\alpha )$ is a summand of $V_{0}^{+}$, i.e. $O(\alpha
)|_{C_{0}}\cong O_{\mathbb{P}^{1}}(1)$, we have $O(\alpha
+C_{0})|_{C_{0}}\cong O_{\mathbb{P}^{1}}(-1)$ and $\alpha +C_{0}=F-(\alpha
_{0}-\alpha )$ with $(\alpha _{0}-\alpha )\in \Phi ^{+}$, that is $O(\alpha
+C_{0})$ is a summand of $V_{0}^{-}(F)$. Since $F=\alpha _{0}+C_{0}$
satisfies $F\cdot F=0$, we have $O(F)|_{C_{0}}\cong O_{\mathbb{P}^{1}}$, $%
O(\alpha _{0})|_{C_{0}}\cong O_{\mathbb{P}^{1}}(2)$ and $O(2F-\alpha
_{0})|_{C_{0}}\cong O_{\mathbb{P}^{1}}(-2)$.

For the loop $E_{8}$-bundle, we have%
\begin{eqnarray*}
\mathcal{E}_{0}^{LE_{8}} &=&\bigoplus_{n\in \mathbb{Z}}(\mathcal{E}%
_{0}^{E_{8}}\otimes O(nF)) \\
&=&\bigoplus_{n\in \mathbb{Z}}(\mathcal{(E}_{0}^{E_{7}}\oplus O\oplus
V_{0}^{+}\oplus V_{0}^{-}\oplus O(\alpha _{0})\oplus O(-\alpha _{0}))\otimes
O(nF)) \\
&=&\bigoplus_{n\in \mathbb{Z}}(\mathcal{(E}_{0}^{E_{7}}\oplus O\oplus
V_{0}^{+}\oplus V_{0}^{-}(F)\oplus O(\alpha _{0}-F)\oplus O(F-\alpha
_{0}))\otimes O(nF))\text{.}
\end{eqnarray*}%
We denote $L_{0}^{248}\triangleq \mathcal{E}_{0}^{E_{7}}\oplus O\oplus
V_{0}^{+}\oplus V_{0}^{-}(F)\oplus O(\alpha _{0}-F)\oplus O(F-\alpha _{0})$.
From definition of $\overline{\partial }_{\varphi }$, $\mathcal{E}_{\varphi
}^{LE_{8}}$ is built from successive extensions of $L_{\varphi
}^{248}\otimes O(nF)$'s, i.e.%
\begin{equation*}
\overline{\partial }_{\varphi }=\left( 
\begin{tabular}{c|c|c|c}
$\ddots $ & $\ddots $ & $\ddots $ & $\ddots $ \\ \hline
$\ddots $ & $\overline{\partial }_{L_{\varphi }^{248}\otimes O((n+1)F)}$ & $%
\ast $ & $\ddots $ \\ \hline
$\ddots $ & $0$ & $\overline{\partial }_{L_{\varphi }^{248}\otimes O(nF)}$ & 
$\ddots $ \\ \hline
$\ddots $ & $\ddots $ & $\ddots $ & $\ddots $%
\end{tabular}%
\ \right) \text{.}
\end{equation*}

So if we can prove $[\varphi _{C_{0}}|_{C_{0}}]\neq 0$ implies $(L_{\varphi
}^{248},\overline{\partial }_{\varphi }|_{L_{\varphi }^{248}})$ is trivial
over $C_{0}$, then $(\mathcal{E}_{\varphi }^{LE_{8}},\overline{\partial }%
_{\varphi })$ is also trivial over $C_{0}$ because of $Ext_{\mathbb{P}%
^{1}}^{1}(O,O)=0$. Note 
\begin{equation*}
L_{0}^{248}|_{C_{0}}\cong O_{\mathbb{P}^{1}}^{\oplus 133}\oplus O_{\mathbb{P}%
^{1}}\oplus (O_{\mathbb{P}^{1}}(1)\oplus O_{\mathbb{P}^{1}}(-1))^{\oplus
56}\oplus O_{\mathbb{P}^{1}}(2)\oplus O_{\mathbb{P}^{1}}(-2)\text{.}
\end{equation*}%
In this decomposition, any of the $56$ pairs of $\{O_{\mathbb{P}^{1}}(-1),O_{%
\mathbb{P}^{1}}(1)\}$ is the restriction of $\{O(\alpha ),O(\alpha
+C_{0})=O(F-(\alpha _{0}-\alpha ))\}$ to $C_{0}$ for some $\alpha $ with $%
a_{1}(\alpha )=1$ and the triple $\{O_{\mathbb{P}^{1}}(2),O_{\mathbb{P}%
^{1}},O_{\mathbb{P}^{1}}(-2)\}$ is the restriction of $%
\{O(-C_{0}),O,O(C_{0})\}$ to $C_{0}$. We will show that the restriction of $%
\overline{\partial }_{\varphi }|_{L_{\varphi }^{248}}$ to $C_{0}$ gives a
non-trivial extension for each of these pairs $\{O_{\mathbb{P}^{1}}(-1),O_{%
\mathbb{P}^{1}}(1)\}$'s and the triple $\{O_{\mathbb{P}^{1}}(-2),O_{\mathbb{P%
}^{1}},O_{\mathbb{P}^{1}}(2)\}$.

In order to write $\overline{\partial }_{\varphi }|_{L_{\varphi }^{248}}$ in
matrix form, we need to decompose $\mathcal{E}_{0}^{E_{7}}$ into positive
parts and non-positive parts, i.e. we denote $\mathcal{E}_{0}^{(E_{7},+)}:=%
\bigoplus_{\alpha \in \Phi ^{+},a_{1}(\alpha )=0}O(\alpha )$ and $\mathcal{E}%
_{0}^{(E_{7},-)}:=O^{\oplus 7}\oplus \bigoplus_{\alpha \in \Phi
^{-},a_{1}(\alpha )=0}O(\alpha )$. Then $\overline{\partial }_{\varphi
}|_{L_{\varphi }^{248}}$ can be written as follows: ($\overline{\partial }%
_{\varphi }|_{L_{\varphi }^{248}}$ is a upper triangle matrix since $%
\overline{\partial }_{\varphi }|_{L_{\varphi }^{248}}$ maps any line bundle
summand to other more "positive" line bundle summands, i.e. $\overline{%
\partial }_{\varphi }:O(D)\rightarrow O(D^{^{\prime }})$ is nonzero only if $%
D^{^{\prime }}-D\geq 0$) 
\begin{equation*}
\overline{\partial }_{\varphi }|_{L_{\varphi }^{248}}=\left( 
\begin{tabular}{c|c|c|c|c|c|c}
$\overline{\partial }_{V_{\varphi }^{-}(F)}$ & $A_{12}$ & $A_{13}$ & $A_{14}$
& $A_{15}$ & $A_{16}$ & $A_{17}$ \\ \hline
$0$ & $\overline{\partial }_{O(F-\alpha _{0})}$ & $A_{23}$ & $A_{24}$ & $%
A_{25}$ & $A_{26}$ & $A_{27}$ \\ \hline
$0$ & $0$ & $\overline{\partial }_{V_{\varphi }^{+}}$ & $A_{34}$ & $A_{35}$
& $A_{36}$ & $A_{37}$ \\ \hline
$0$ & $0$ & $0$ & $\overline{\partial }_{\mathcal{E}_{\varphi }^{(E_{7},+)}}$
& $A_{45}$ & $A_{46}$ & $A_{47}$ \\ \hline
$0$ & $0$ & $0$ & $0$ & $\overline{\partial }_{O}$ & $A_{56}$ & $A_{57}$ \\ 
\hline
$0$ & $0$ & $0$ & $0$ & $0$ & $\overline{\partial }_{\mathcal{E}_{\varphi
}^{(E_{7},-)}}$ & $A_{67}$ \\ \hline
$0$ & $0$ & $0$ & $0$ & $0$ & $0$ & $\overline{\partial }_{O(\alpha _{0}-F)}$%
\end{tabular}%
\ \right) \text{.}
\end{equation*}

Now we restrict this to $C_{0}$, the $56$ pairs $\{O_{\mathbb{P}^{1}}(-1),O_{%
\mathbb{P}^{1}}(1)\}$'s are in $V_{0}^{-}(F)|_{C_{0}}\oplus
V_{0}^{+}|_{C_{0}}$. Since $A_{23}=(0,0,\cdots ,0)_{56\times 1}$ and%
\begin{equation*}
A_{13}=\left( 
\begin{array}{cccc}
\pm \varphi _{C_{0}} & \ast & \cdots & \ast \\ 
0 & \pm \varphi _{C_{0}} & \cdots & \ast \\ 
\vdots & \vdots & \ddots & \vdots \\ 
0 & 0 & \cdots & \pm \varphi _{C_{0}}%
\end{array}%
\right) _{56\times 56}\text{,}
\end{equation*}%
if $[\varphi _{C_{0}}|_{C_{0}}]\neq 0$, then we have a trivialization of the 
$56$ pairs $\{O_{\mathbb{P}^{1}}(-1),O_{\mathbb{P}^{1}}(1)\}$'s over $C_{0}$
by Lemma $32$ in \cite{CL}.

For the triple $\{O_{\mathbb{P}^{1}}(-2),O_{\mathbb{P}^{1}},O_{\mathbb{P}%
^{1}}(2)\}$, we review the trivialization of $A_{1}$ Lie algebra bundle. In $%
A_{1}$ case, we have an $A_{1}$-bundle $\mathcal{E}_{\varphi }^{A_{1}}$,
which topologically is $\mathcal{E}_{0}^{A_{1}}=O\oplus O(C)\oplus O(-C)$,
but with a holomorphic structure as follows:%
\begin{equation*}
\overline{\partial }_{\varphi }=\left( 
\begin{tabular}{c|c|c}
$\overline{\partial }_{0}$ & $\pm \varphi _{C}$ & $0$ \\ \hline
$0$ & $\overline{\partial }_{0}$ & $\pm \varphi _{C}$ \\ \hline
$0$ & $0$ & $\overline{\partial }_{0}$%
\end{tabular}%
\ \right) \text{,}
\end{equation*}%
where $\varphi _{C}\in H^{0,1}(X,O(C))$. From \cite{CL}, we know if $%
[\varphi _{C}|_{C}]\neq 0$, then $\mathcal{E}_{\varphi }^{A_{1}}$ is trivial
on $C$. Back to our case, the triple $\{O_{\mathbb{P}^{1}}(-2),O_{\mathbb{P}%
^{1}},O_{\mathbb{P}^{1}}(2)\}$ has the corresponding submatrices $%
A_{25}=(\varphi _{C_{0}})_{1\times 1}$, $A_{57}=(\varphi _{C_{0}})_{1\times
1}$ and $A_{27}=(0)_{1\times 1}$. Since $A_{23}$, $A_{24}$, $A_{26}$, $%
A_{47} $ and $A_{67}$ are all zero matrices, from the trivialization of $%
A_{1}$ Lie algebra bundle, we know if $[\varphi _{C_{0}}|_{C_{0}}]\neq 0$,
then we have a trivialization of the triple $\{O_{\mathbb{P}^{1}}(-2),O_{%
\mathbb{P}^{1}},O_{\mathbb{P}^{1}}(2)\}$ over $C_{0}$.

Hence if $[\varphi _{C_{0}}|_{C_{0}}]\neq 0$, then $(L_{\varphi }^{248},%
\overline{\partial }_{\varphi }|_{L_{\varphi }^{248}})$ is trivial on $C_{0}$%
, which implies $(\mathcal{E}_{\varphi }^{LE_{8}},\overline{\partial }%
_{\varphi })$ is also trivial on $C_{0}$. Hence, we have Theorem \ref{loop
thm} for $LE_{8}$ case.

\section{$E_{n}$-bundle over $X_{n}$ with $n\leq 9$}

When $X=X_{n}$ is a blowup of $\mathbb{P}^{2}$ at $n$ points $x_{1},\cdots
,x_{n}$ with $n\leq 9$, there is a canonical (affine) Lie algebra bundle $%
\mathcal{E}_{0}^{E_{n}}$ over it, where $E_{9}$ is the affine $E_{8}$. In
this section, we will give a detail study of the relationship between the
geometry of $X_{n}$ and the deformability of $\mathcal{E}_{0}^{E_{n}}$.

\subsection{$E_{n}$-bundle over $X_{n}$ with $n\leq 9$}

The Picard group $Pic(X_{n})\cong H^{2}(X_{n},\mathbb{Z)}$ is a rank $n+1$
lattice with generators $h,l_{1},\cdots ,l_{n}$, where $h$ is the class of
lines in $\mathbb{P}^{2}$ and $l_{i}$ is the exceptional class of the
blow-up at $x_{i}$. So $h^{2}=1=-l_{i}^{2}$ and $h\cdot l_{i}=0=l_{i}\cdot
l_{j}$, $i\neq j$. Thus $H^{2}(X_{n},\mathbb{Z)\cong Z}^{1,n}$. The
canonical class is $K_{X_{n}}=-3h+l_{1}+\cdots +l_{n}$. Denote%
\begin{equation*}
\Phi _{n}:=\{\alpha \in H^{2}(X_{n},\mathbb{Z)}|\alpha ^{2}=-2,\alpha \cdot
K=0\}\text{.}
\end{equation*}%
Then $\Phi _{n}$ is a root system of type $E_{n}$ when $n\leq 8$ and $\Phi
_{9}$ is an affine real root system of $\widehat{E}_{8}$ (also denoted as $%
E_{9}$). More explicitly, $\Phi _{\widehat{E}_{8}}:=\Phi _{9}\cup
\{mK_{X_{9}}|m\neq 0,m\in \mathbb{Z}\}$ forms a root system of (untwisted)
affine $E_{8}$-type (that is, $\widehat{E}_{8}$-type) with $\Phi _{\widehat{E%
}_{8}}^{re}:=\Phi _{9}$ the set of real roots and $\Phi _{\widehat{E}%
_{8}}^{im}:=\{mK_{X_{9}}|m\neq 0,m\in \mathbb{Z}\}$ the set of imaginary
roots (see \cite{HL} or \cite{LXZ}). We have an $\widehat{E}_{8}$-bundle $%
\mathcal{E}_{0}^{\widehat{E}_{8}}$ over $X_{9}$:%
\begin{equation*}
\mathcal{E}_{0}^{\widehat{E}_{8}}=O^{\oplus 9}\oplus \bigoplus_{\alpha \in
\Phi _{\widehat{E}_{8}}^{re}}O(\alpha )\bigoplus_{\beta \in \Phi _{\widehat{E%
}_{8}}^{im}}O(\beta )
\end{equation*}%
The Lie algebra structure on $\mathcal{E}_{0}^{\widehat{E}_{8}}$ is
explained in \cite{LXZ}. When $n\leq 8$, $\mathcal{E}_{0}^{E_{n}}=O^{\oplus
n}\oplus \bigoplus_{\alpha \in \Phi _{n}}O(\alpha )$ is an $E_{n}$-bundle
over $X_{n}$.

Suppose $C=\cup C_{i}$ is an (affine) $ADE$ curve of type $\mathfrak{g}$ in $%
X_{n}$, then $C_{i}$'s generates a subroot system $\Phi $ inside $\Phi _{n}$
since $C_{i}\cdot K=0$ for every $i$. Therefore the corresponding bundle $%
\mathcal{E}_{0}^{\mathfrak{g}}$ is a Lie algebra subbundle of $\mathcal{E}%
_{0}^{E_{n}}$.

Suppose $\mathcal{E}_{0}^{\mathfrak{g}}$ is a $\mathfrak{g}$-bundle over a
surface $X$ corresponding to a root system $\Lambda _{\mathfrak{g}}\subset
Pic(X)$ of type $\mathfrak{g}$.

\begin{definition}
A Lie algebra sub-bundle $\mathcal{F}$ of $\mathcal{E}_{0}^{\mathfrak{g}}$
is called strict if there exists a sub-root lattice $\Lambda $ of $\Lambda _{%
\mathfrak{g}}$ such that $\mathcal{F}$ is a direct sum of line bundles
corresponding to the roots in $\Lambda $.
\end{definition}

In order to describe $\mathcal{E}_{0}^{\widehat{E}_{8}}$ as a central
extension of a loop Lie algebra bundle over $X_{9}$, we pick any smooth $%
\left( -1\right) $-curve $l$ in $X_{9}$, then we have%
\begin{equation*}
\mathcal{E}_{0}^{\widehat{E}_{8}}\cong \mathcal{E}_{0}^{E_{8}}\otimes
(\bigoplus_{n\in \mathbb{Z}}O(nK_{X_{9}}))\oplus O\text{,}
\end{equation*}%
where $\mathcal{E}_{0}^{E_{8}}$ is the pull-back of the $E_{8}$-bundle over $%
X_{8}$ via $\pi :X_{9}\rightarrow X_{8}$, the blow down map of $l$. The next
proposition describes the converse.

\begin{proposition}
When $\mathcal{E}_{0}^{\widehat{E}_{8}}$ is a central extension of a loop $%
E_{8}$-sub-bundle over $X$ for some strict $E_{8}$-bundle $\mathcal{F}%
_{0}^{E_{8}}$ over $X_{9}$, i.e.%
\begin{equation*}
\mathcal{E}_{0}^{\widehat{E}_{8}}\cong \mathcal{F}_{0}^{E_{8}}\otimes
(\bigoplus_{n\in \mathbb{Z}}O(nK_{X_{9}}))\oplus O\text{,}
\end{equation*}%
as a Lie algebra bundle isomorphism, then there is a unique $($possibly
reducible$)$ $\left( -1\right) $-curve $l$ in $X$ such that $\mathcal{F}%
_{0}^{E_{8}}$ is constructed from those $\alpha \in \Lambda ^{re}$
satisfying $\alpha \cdot l=0$.

\begin{proof}
Denote $\Delta _{E_{8}}=\{\alpha _{1},\cdots ,\alpha _{8}\}$ as a root base
of the corresponding $E_{8}$ Lie algebra from $\mathcal{F}_{0}^{E_{8}}$, we
need to find a unique $(-1)$-curve $l$ in $X$ such that $l\cdot \alpha
_{i}=0 $ for any $\alpha _{i}$ in $\Delta _{E_{8}}$. Since $\{\pm 1\}\times
W(\widehat{E}_{8})$ acts on the set of all root bases of $\widehat{E}_{8}$
simply transitively \cite{Kac} and $W(\widehat{E}_{8})$ acts on the set of $%
(-1)$-curves \cite{LXZ}, we only need to find $l$ for one particular root
base of any $E_{8}$ in $\widehat{E}_{8}$ and show that such a $l$ is unique.
For example, if we take $\alpha _{1}=h-l_{1}-l_{2}-l_{3},\alpha
_{k}=l_{k-1}-l_{k}$ for $k=2,\cdots 8$, then we can take $l=l_{9}$ and by
the condition that $l\cdot \alpha _{i}=0$, $l^{2}=-1=l\cdot K$, we know such
a $l$ is unique.
\end{proof}
\end{proposition}

\subsection{Deformability of such $\mathcal{E}_{0}^{\widehat{E}_{8}}$}

In this subsection, we will describe relationships between the geometry of $%
X_{9}$ and the deformability of $\mathcal{E}_{0}^{\widehat{E}_{8}}$. Similar
results for $X_{n}$ and $\mathcal{E}_{0}^{E_{n}}$ with $n\leq 8$ can be
easily deduced from this case.

Recall when $Pic(X)$ contains a lattice $\Lambda $ isomorphic to a root
lattice $\Lambda _{\mathfrak{g}}$, then we have a $\mathfrak{g}$-bundle $%
\mathcal{E}$ over $X$ (\cite{D}\cite{FMW}\cite{LZ}\cite{LZ2}\cite{LXZ}).%
\begin{equation*}
\mathcal{E}:=O^{\oplus r}\oplus \bigoplus_{\alpha \in \Phi }O(\alpha )\text{.%
}
\end{equation*}%
Infinitesimal deformations of holomorphic structures on $\mathcal{E}$ are
parametrized by $H^{1}(X,End(\mathcal{E}))$, and those which also preserve
the Lie algebra structure are parametrized by $H^{1}(X,ad(\mathcal{E}%
))=H^{1}(X,\mathcal{E})$ since $\mathfrak{g}$ is simple. Hence we introduce
the following definitions.

\begin{definition}
\label{deform}$(i)$ $\mathcal{E}$ is called fully deformable if there exists
a base $\Delta \subset \Phi $ such that $H^{1}(X,O(\alpha ))\neq 0$ for any $%
\alpha \in \Delta $.

$(ii)$ $\mathcal{E}$ is called $\mathfrak{h}$-deformable if there exists a
strict $\mathfrak{h}$ Lie algebra sub-bundle $\mathcal{E}^{\mathfrak{h}%
}\subseteq \mathcal{E}$ which is fully deformable.

$(iii)$ $\mathcal{E}$ is called deformable in $\alpha $-direction for $%
\alpha \in \Phi $ if $H^{1}(X,O(\alpha ))\neq 0$.

$(iv)$ $\mathcal{E}$ is called totally non-deformable if $H^{1}(X,O(\alpha
))=0$ for any $\alpha \in \Phi $.
\end{definition}

Recall the holomorphic structure $\overline{\partial }_{\varphi }$ or $%
\overline{\partial }_{(\varphi ,\Phi )}$ defined in $\S 3.1$ and $\S 3.2$ on 
$\mathcal{E}$ admits a filtration determined by the height of the roots (if
the root base $\Delta =\{\alpha _{1},\alpha _{2},\cdots ,\alpha _{r}\}$,
then for any $\alpha \in \Phi $, we have $\alpha =\sum a_{i}\alpha _{i}$ and
the height of $\alpha $ is defined to be $ht(\alpha ):=\sum a_{i}$).

\begin{remark}
When $\mathcal{E}$ is fully deformable and if for every simple root $\alpha
\in \Delta $, $O(\alpha )=O(C_{\alpha })$ for some smooth irreducible curve $%
C_{\alpha }$, then $C=\cup _{\alpha \in \Delta }C_{\alpha }$ is an $ADE$ or
affine $ADE$ curve in $X$. In this case, we can show that $H^{2}(X,O(\alpha
))=0$ for any $\alpha \in \Phi $ and the $\mathfrak{g}$ or $\widehat{%
\mathfrak{g}}$ bundle $\mathcal{E}$ admits a deformation into a filtrated
bundle which is trivial on every $C_{\alpha }$ $($see section $3)$. When $%
\mathcal{E}$ is totally non-deformable, $\overline{\partial }_{\varphi }$
can only be $\overline{\partial }_{0}$.
\end{remark}

The main results of this section are the followings.

\begin{theorem}
\label{general position}$\mathcal{E}_{0}^{\widehat{E}_{8}}$ over $X_{9}$ is
totally non-deformable if and only if the nine blowup points in $\mathbb{P}%
^{2}$ are in general position.
\end{theorem}

Let us recall some facts about elliptic fibrations on $X_{9}$ \cite{M}\cite%
{P}. Any elliptic fibration on $X_{9}$ must be relatively minimal, i.e.
there is no $(-1)$-curves in any of its fibrations, as there is no elliptic
fibration on $X_{8}$, this is because the Euler characteristic of any
elliptic surface is a multiple of $12$ \cite{F} and also $\chi (X_{9})=12$.
There is at most one multiple fiber \cite{F2}, say of multiplicity $m$. This
happens precisely when there exists an irreducible pencil of degree $3m$ in $%
\mathbb{P}^{2}$ with $9$ base points, each of multiplicity $m$ and $X_{9}$
is the blow up of $\mathbb{P}^{2}$ at these $9$ points. We can characterize
the existence of such an elliptic fibration on $X_{9}$ in terms of
deformability of $\mathcal{E}_{0}^{\widehat{E}_{8}}$ along imaginary root
directions. For instance, $X_{9}$ with $-K_{X_{9}}$ nef admits an elliptic
fibration (without multiple fiber) if and only if $\mathcal{E}_{0}^{\widehat{%
E}_{8}}$ is deformable in $(-mK)$-direction for some $m\in \mathbb{N}$ (with 
$m=1$). Deformability of $\mathcal{E}_{0}^{\widehat{E}_{8}}$ can also detect
the existence of $ADE$ or Kodaira curves in $X$.

\begin{theorem}
\label{elliptic}Suppose $-K_{X_{9}}$ is nef, then

$(i)$ $X_{9}$ admits an elliptic fibration with a multiple fiber of
multiplicity $m$ $\left( m\geq 1\right) $ if and only if $\mathcal{E}_{0}^{%
\widehat{E}_{8}}$ is deformable in $(-mK)$-direction but not in $(-m+1)K$%
-direction.

$(ii)$ $X_{9}$ has an $($maximal$)$ $ADE$ curve $C$ of type $\mathfrak{g}$
if and only if $\mathcal{E}_{0}^{\widehat{E}_{8}}$ is $($maximal$)$ $%
\mathfrak{g}$-deformable.

$(iii)$ $X_{9}$ has a $($maximal$)$ Kodaira curve $C$ of type $\widehat{%
\mathfrak{g}}$ if and only if $\mathcal{E}_{0}^{\widehat{E}_{8}}$ is $($%
maximal$)$ $\widehat{\mathfrak{g}}$-deformable.
\end{theorem}

Here we say an $ADE$ or Kodaira curve $C$ is maximal if it is not proper
contained in another $ADE$ or Kodaira curve. We say $\mathcal{E}_{0}^{%
\widehat{E}_{8}}$ is maximal $\mathfrak{g}$ (or $\widehat{\mathfrak{g}}$)
deformable if there does not exist another fully deformable $($affine$)$ Lie
algebra sub-bundle of $\mathcal{E}_{0}^{\widehat{E}_{8}}$ containing this $%
\mathfrak{g}$ (or $\widehat{\mathfrak{g}}$) bundle.

\subsection{Negative curves in $X_{9}$}

In this subsection, we study negative rational curves in $X_{9}$. We can get
corresponding results for $X_{n}$ with $n\leq 8$ from this $n=9$ case.

A divisor $D$ in $X$ is called a $\left( -m\right) $-class if $D\cdot D=-m$
and $D\cdot K=m-2$. An effective $\left( -m\right) $-class is called a $(-m)$%
-curve. Note when $D=\sum n_{i}C_{i}$ is a $(-m)$-curve, we will also denote
the corresponding curve $\cup C_{i}$ as $D$.

Use the notations in the above subsection, every effective divisor $%
D=ah-\sum_{i=1}^{9}a_{i}l_{i}\in Pic(X_{9})$ must have $a=D\cdot h\geq 0$.
It is well-known that all $(-1)$-classes are effective, and there are
infinite number of them in $X_{9}$. There are also infinite number of $(-2)$%
-classes, but whether they are effective or not depends on the positions of
the $9$ blow-up points.

\begin{definition}
Let $x_{1},\cdots ,x_{n}$ be $n$ distinct points in $\mathbb{P}^{2}$. These $%
n$ points are said to be non-special with respect to Cremona transformations
if for any Cremona transformation $T$ with centers within $x_{i}$'s, the
points $y_{1},\cdots ,y_{n}$ corresponding to $x_{i}$'s under $T$ are
distinct points such that no three points among $y_{1},\cdots ,y_{n}$ are
collinear.
\end{definition}

\begin{definition}
$($\cite{LXZ}$)$ Let $x_{1},\cdots ,x_{9}$ be $9$ points in $\mathbb{P}^{2}$%
, we say they are in general position if they satisfy the following three
conditions:

$(i)$ they are distinct points in $\mathbb{P}^{2}$;

$(ii)$ they are non-special with respect to Cremona transformations;

$(iii)$ there is a unique cubic curve passing through all of them.
\end{definition}

The conditions (i) and (ii) mean that any $8$ of these $9$ points are in
general position. That is, no lines pass through three of them, no conics
pass through six of them, and no cubic curves pass through eight of them
with one of the eight points being a double point.

If the $9$ blowing up points are in general position, then there is no
effective $(-2)$-class in $X_{9}$ \cite{LXZ}. In general, there are at most
finite number of $(-m)$-curves with $m\geq 3$.

\begin{lemma}
\label{curve}Let $D=ah-\sum_{i=1}^{9}a_{i}l_{i}$ be a $(-m)$-curve in $X_{9}$
with $m\geq 3$, then

$(i)$ $m\leq 9$;

$(ii)$ $0\leq a\leq 3$;

$(iii)$ $-1\leq a_{i}\leq 2$ for all $i$, and there exists some $j$ with $%
a_{j}=1$;

$(iv)$ there are finite number of such curves.

\begin{proof}
$(i)$ Since $D$ is a $(-m)$-curve, $D\cdot D=-m$ and $D\cdot K=m-2$, i.e.%
\begin{equation*}
\sum a_{i}^{2}=a^{2}+m\text{ and }\sum a_{i}=3a+m-2\text{.}
\end{equation*}%
From the above two equations, we have%
\begin{equation*}
(3a+m-2)^{2}=(\sum a_{i})^{2}\leq 9(\sum a_{i}^{2})=9(a^{2}+m)\text{.}
\end{equation*}%
Thus, $a\leq \frac{-m^{2}+13m-4}{6(m-2)}$, also $a\geq 0$ since $D$ is
effective, hence $m\leq 12$.

When $m\geq 10$, we must have $a=0$, that means $\sum a_{i}^{2}=m$ and $\sum
a_{i}=m-2$, hence $\sum a_{i}^{2}-\sum a_{i}=2$, which implies every $a_{i}$
satisfies $|a_{i}|\leq 1$ and there exists exactly one $a_{i}$ with $%
a_{i}=-1 $. But we also have $\sum a_{i}=m-2\geq 8$, which is impossible
since we only have nine $a_{i}$'s.

$(ii)$ When $m\geq 4$, $a\leq \frac{-m^{2}+13m-4}{6(m-2)}\leq \frac{8}{3}<3$%
. When $m=3$, $a\leq \frac{-m^{2}+13m-4}{6(m-2)}=\frac{13}{3}<5$. Hence we
only need to prove there is no $(-3)$-curve with $a=4$.

Suppose not, then there exists $a_{i}$'s such that $\sum a_{i}^{2}=19$ and $%
\sum a_{i}=13$. From $\sum a_{i}^{2}-\sum a_{i}=6$, we know $-2\leq
a_{i}\leq 3$. If there is any $a_{i}$ with $a_{i}=3$, then the other $a_{i}$%
's can only be $0$ or $1$, but we have $\sum a_{i}=13$ and there is only
nine $a_{i}$'s, which is impossible. Hence $-2\leq a_{i}\leq 2$, from $\sum
a_{i}^{2}-\sum a_{i}=6$, we can have at most three $a_{i}$'s equal to $2$,
which is also impossible since $\sum a_{i}=13$.

$(iii)$ From $\sum a_{i}^{2}=a^{2}+m$, $\sum a_{i}=3a+m-2$ and $0\leq a\leq
3 $, we have%
\begin{equation*}
\sum a_{i}=3a+m-2\geq a^{2}+m-2=\sum a_{i}^{2}-2\text{.}
\end{equation*}%
Hence $-1\leq a_{i}\leq 2$. And there are three cases:

Case 1, one $a_{i}$ equal to $2$, the others equal to $0$ or $1$;

Case 2, one $a_{i}$ equal to $-1$, the others equal to $0$ or $1$;

Case 3, all $a_{i}$'s are equal to $0$ or $1$.

By $\sum a_{i}=3a+m-2\geq 1$, we know in case 2 and case 3, there must exist
some $a_{i}$ with $a_{i}=1$. In case 1, if there is no $a_{i}$ with $a_{i}=1$%
, then $D=ah-2l_{j}$. From $\sum a_{i}^{2}=a^{2}+m$, $\sum a_{i}=3a+m-2$, we
have $a=0$, $m=4$, hence $D=-2l_{j}$, which is not an effective divisor.

$(iv)$ It is obvious from the above results.
\end{proof}
\end{lemma}

From this lemma, we can easily obtain the following as a corollary.

\begin{corollary}
If there exists a $(-m)$-curve in $X_{9}$ with $m\geq 3$, then there also
exists a $(-m+1)$-curve in $X_{9}$.

\begin{proof}
If $D\in |ah-\sum a_{i}l_{i}|$ is a $(-m)$-curve in $X_{9}$ with $m\geq 3$,
then there exists $j$ with $a_{j}=1$ by $(iii)$ of Lemma \textit{\ref{curve}}%
. It is easy to check that $D+l_{j}$ is a $(-m+1)$-curve in $X_{9}$.
\end{proof}
\end{corollary}

If the $9$ blowing up points are in general position, then there is no $(-2)$%
-curve in $X_{9}$, as a consequence, there is also no $(-m)$-curve in $X_{9}$
with $m\geq 3$. The following result shows that this happens exactly when $%
X_{9}$ is almost Fano. We include a proof here as we could not find it in
the literatures.

\begin{lemma}
\label{nef}$X_{9}$ has no $(-m)$-curve with $m\geq 3$ if and only if $%
-K_{X_{9}}$ is nef.

\begin{proof}
If $-K$ is nef, then from $C\cdot K^{-1}=2-m\geq 0$ for any $(-m)$-curve $C$%
, we know $m\leq 2$.

Conversely, assume $X_{9}$ has no $(-m)$-curve with $m\geq 3$. Since $X_{9}$
is a blowup of $\mathbb{P}^{2}$ at nine points $\{x_{i}\}_{i=1}^{9}$, we
have an effective anti-canonical divisor $D$. Recall when $D\cdot \Sigma <0$
for any irreducible curve $\Sigma $ in $X$, $\Sigma $ must be a component of 
$D$. So if $D$ is an irreducible curve or a Kodaira curve, then $D$ is nef.
We denote the image of $D$ in $\mathbb{P}^{2}$ as $C$, which is a cubic
curve passing through these $9$ blowing up points.

$(i)$ If $C$ is smooth, then we are done as $D\cong C$ and therefore
irreducible.

$(ii)$ If $C$ is reduced and irreducible, then it must be a nodal or
cuspidal cubic. If $\{x_{i}\}_{i=1}^{9}\cap $sing$(C)=\varnothing $ $%
(sing(C) $ means the set of singular points on $C)$, then $D\cong C$ and we
are done. Otherwise, say $x_{1}\in $sing$(C)$ and we write the strict and
proper transformations of $C$ in $Bl_{x_{1}}(\mathbb{P}^{2})$ as $C_{1}$ and 
$C_{1}+E$ respectively. Then the remaining $x_{i}$'s must have exactly $1$
point $($resp. $7$ points$)$ lying on $E$ $($resp. $C_{1})$ in order to
avoid having $(-m)$-curve with $m\geq 3$. Thus $D$ is a Kodaira curve of
type $\widehat{A}_{1}$ or $III(\widehat{A}_{1})$ for $C$ being a nodal or
cuspidal respectively.

$(iii)$ If $C$ is reduced and reducible, then $C=B\cup H_{0}$ or $H_{1}\cup
H_{2}\cup H_{3}$ with $B$ and $H_{j}$'s are conic and distinct lines in $%
\mathbb{P}^{2}$. As before, we must have exactly $6$ $x_{i}$'s on $B$ and $3$
$x_{i}$'s on each $H_{j}$ and none on sing$(C)$. Thus $D\cong C$ is a
Kodaira curve of type $\widehat{A}_{1}$, $\widehat{A}_{2}$, $III(\widehat{A}%
_{1})$ or $VI(\widehat{A}_{2})$.

$(iv)$ If $C$ is non-reduced, $C=3H$, $D$ must have a $(-m)$-curve with $%
m\geq 3$.

Hence $D$ is an irreducible curve or a Kodaira curve, and we are done.
\end{proof}
\end{lemma}

In the following two lemmas, we will use Lemma $2.21$ in \cite{BPV} to give
a criteria of a curve in $X_{n}$ being an $ADE$ or affine $ADE$ curve. Lemma 
$2.21$ can be reformulated as follows: if $C=\cup _{i=1}^{r}C_{i}$ is a
connected curve in a surface $X$ satisfying: (i) $C_{i}^{2}=-2$ and $%
C_{i}\cdot K_{X}=0$ for any $i$; (ii) $C_{i}\cdot C_{j}\leq 1$ for any $%
i\neq j$; (iii) $(C_{i}\cdot C_{j})_{r\times r}\leq 0$. Then when $%
(C_{i}\cdot C_{j})_{r\times r}<0$, $C$ is an $ADE$ curve, otherwise, it is
an affine $ADE$ curve.

\begin{lemma}
Suppose $-K_{X_{n}}$ $(n\leq 8)$ is nef. Let $C=\cup C_{i}$ be a connected
curve in $X_{n}$. If $C\cdot K_{X_{n}}=0$, then $C$ is an $ADE$ curve.

\begin{proof}
Since $-K_{X_{n}}$ is nef, $C\cdot K_{X_{n}}=0$ implies $C_{i}\cdot
K_{X_{n}}=0$ for each $i$, i.e. $[C_{i}]\in \langle K\rangle ^{\perp }\cong
\Lambda _{E_{n}}$. We have $C_{i}^{2}<0$ and $(C_{i}+C_{j})^{2}<0$ for any $%
i $ and $j$. Together with the genus formula, we have $C_{i}^{2}=-2$ and $%
C_{i}\cdot C_{j}\leq 1$ for $i\neq j$. By Lemma $2.21$ in \cite{BPV}\textit{%
, we know }$C$ is an $ADE$ curve.
\end{proof}
\end{lemma}

For $n=9$ case, we have the following lemma.

\begin{lemma}
\label{curve2}Suppose $-K_{X_{9}}$ is nef. Let $C=\cup C_{i}$ be a connected
curve in $X_{9}$. If $C\cdot K_{X_{9}}=0$ and $C_{i}+K_{X_{9}}$ is not
effective for each $i$, then $C$ is a smooth elliptic curve, an $ADE$ curve
or an affine $ADE$ curve.

\begin{proof}
Since $-K_{X_{9}}$ is nef, $C\cdot K_{X_{9}}=0$ implies $C_{i}\cdot
K_{X_{9}}=0$ for each $i$, i.e. $[C_{i}]\in \langle K_{X_{9}}\rangle ^{\perp
}\cong \Lambda _{E_{9}}$. We have $C_{i}^{2}\leq 0$ and $(C_{i}+C_{j})^{2}%
\leq 0$ for any $i$ and $j$. Moreover, for any effective divisor $D\in
\langle K_{X_{9}}\rangle ^{\perp }$, if $D^{2}=0$, then $D\in |mK_{X_{9}}|$
for some non-zero integer $m$. From $C_{i}^{2}\leq 0$ and genus formula, we
have $C_{i}^{2}=-2$ or $0$.

If there exists $C_{i}$ such that $C_{i}^{2}=0$, then $C_{i}\in |mK|$ for
some non-zero integer $m$. Since $C_{i}+K_{X_{9}}$ is not effective, we know 
$m=-1$, i.e. $C_{i}\in |-K|$. If $C$ is not irreducible, then there exists $%
C_{j}$ which intersects $C_{i}$, which is impossible. So $C=C_{i}\in |-K|$
is an elliptic curve or an affine $A_{0}$ curve by Lemma \ref{nef}.

If $C_{i}^{2}=-2$ for any $i$, then $C_{i}\cdot C_{j}\leq 2$ for any $i\neq
j $. If there exist $C_{i}$ and $C_{j}$ such that $C_{i}\cdot C_{j}=2$, then 
$(C_{i}+C_{j})^{2}=0$, $C_{i}+C_{j}\in |mK|$ for some integer $m$. Hence $%
C=C_{i}\cup C_{j}$ is an affine $A_{1}$ curve, this is because if $C_{k}$ is
another irreducible component of $C$ and assume it intersects with $C_{i}$,
then it must be an irreducible component of $C_{j}$, which contradicts to $%
C_{j}$ being irreducible. Otherwise, we will have $C_{i}^{2}=-2$ for each $i$
and $C_{i}\cdot C_{j}\leq 1$ for $i\neq j$. By Lemma $2.21$ of \textit{\cite%
{BPV}, we know }$C$ is an $ADE$ or affine $ADE$ curve.
\end{proof}
\end{lemma}

\subsection{Proof of theorems \textbf{\protect\ref{general position} and 
\protect\ref{elliptic}}}

\begin{proof}
$\mathit{(}$\textit{of Theorem \ref{general position}}$\mathit{)}$\textit{\
If the nine blowup points in }$\mathbb{P}^{2}$\textit{\ are in general
position, then for any }$\alpha \in \Phi _{9}$\textit{, we have }$%
h^{0}\left( X,O\left( \alpha \right) \right) =0$\textit{\ \cite{LXZ}. Since }%
$K\cdot K=0$\textit{, we also have }$K-\alpha \in \Phi _{9}$\textit{\ and
therefore }$h^{2}\left( X,O\left( \alpha \right) \right) =0$\textit{\ by
Serre duality. However the Riemann-Roch formula gives }$\chi \left(
X,O\left( \alpha \right) \right) =1+\frac{\alpha ^{2}-\alpha K}{2}=0$\textit{%
\ and therefore }$h^{1}\left( X,O\left( \alpha \right) \right) =0$\textit{.
For the imaginary roots }$mK$\textit{'s, from Lemma }$4$\textit{\ and
Proposition }$11$\textit{\ in \cite{LXZ}, we have }$h^{0}(X,O(mK))=0$\textit{%
\ and }$h^{0}(X,O(-mK))=1$\textit{\ for }$m\geq 1$\textit{. By Serre duality
and Riemann-Roch formula, we have }$h^{1}(X,O(mK))=0$\textit{\ for any
imaginary root }$mK$\textit{. Hence }$\mathcal{E}_{0}^{\widehat{E}_{8}}$%
\textit{\ is totally non-deformable.}

\textit{Conversely, if }$\mathcal{E}_{0}^{\widehat{E}_{8}}$\textit{\ is
totally non-deformable, then }$X$\textit{\ has no }$\mathit{(}$\textit{%
possibly reducible}$\mathit{)}$\textit{\ }$(-2)$\textit{-curve, hence no }$%
(-n)$\textit{-curve with }$n\geq 2$\textit{. By Proposition }$10$\textit{\
in \cite{N}, this implies the nine blowup points are non-special with
respect to Cremona transformations. Also from }$h^{1}(X,O(mK))=0$\textit{\
for any imaginary root }$mK$\textit{, we get }$h^{0}(X,O(-K))=1$, \textit{we
have a unique cubic curve in }$\mathbb{P}^{2}$\textit{\ passing through all
of the blow-up points. Hence, the nine blow-up points in }$\mathbb{P}^{2}$%
\textit{\ are in general position.}
\end{proof}

\bigskip

\begin{proof}
$\mathit{(}$\textit{of Theorem \ref{elliptic}}$\mathit{)}$\textit{\ }$%
\mathit{(i)}$\textit{We have }$h^{1}(X,O(-mK))=h^{0}(X,O(-mK))-1$\textit{\
for any }$m$ \textit{by Riemann-Roch formula. So }$\mathcal{E}_{0}^{\widehat{%
E}_{8}}$\textit{\ is deformable in }$(-mK)$\textit{-direction if and only if 
}$h^{0}(X,O(-mK))=2$.

\textit{Let }$F_{0}\in |-K|$\textit{, then by Proposition }$2.2$\textit{\ of 
\cite{CD},\ }$X$\textit{\ admits an elliptic fibration with a multiple fiber
of multiplicity }$m$ \textit{if and only if }$O_{F_{0}}(F_{0})$ \textit{is
of order }$m$ \textit{in }$Pic(F_{0})$\textit{. But }$O_{F_{0}}(mF_{0})\cong
O_{F_{0}}$ \textit{if and only if }$h^{0}(O_{F_{0}}(mF_{0}))=1$ \textit{as }$%
O_{F_{0}}(mF_{0})$ \textit{is topologically trivial. By the exact sequence}%
\begin{equation*}
0\longrightarrow O_{X}\longrightarrow O_{X}(mF_{0})\longrightarrow
O_{F_{0}}(mF_{0})\longrightarrow 0
\end{equation*}%
\textit{together with }$h^{1}(X,O_{X})=0$\textit{, we know }$%
h^{0}(O_{X}(mF_{0}))=1+h^{0}(O_{F_{0}}(mF_{0}))$\textit{. So } $m=\min
\{n:h^{0}(O_{F_{0}}(nF_{0}))=1\}=\min \{n:h^{0}(X,O(-nK))=2\}$.

$(ii)$\textit{\ If }$X$\textit{\ has an }$ADE$\textit{\ curve }$C$\textit{\
of type }$\mathfrak{g}$\textit{, we can use it to construct a fully
deformable }$\mathfrak{g}$-\textit{subbundle of }$\mathcal{E}_{0}^{\widehat{E%
}_{8}}$\textit{\ as in }$\S 3.2$\textit{. When }$C$\textit{\ is maximal,
then this }$\mathfrak{g}$-\textit{subbundle is not contained in any other
fully deformable Lie algebra subbundle of }$\mathcal{E}_{0}^{\widehat{E}%
_{8}} $.

\textit{Conversely, if }$\mathcal{E}_{0}^{\widehat{E}_{8}}$ \textit{is
maximal }$\mathfrak{g}$\textit{-deformable, then we can find a base }$\Delta
\subset $ $\Phi _{\widehat{E}_{8}}$\textit{\ of }$\mathfrak{g}$\textit{\
such that }$h^{1}(X,O(\alpha ))\neq 0$\textit{\ for every }$\alpha \in
\Delta $\textit{. Since }$\chi (O(\alpha ))=1+\frac{\alpha ^{2}-\alpha \cdot
K}{2}=0$\textit{, we must have }$h^{0}(O(\alpha ))\neq 0$\textit{\ or }$%
h^{2}(O(\alpha ))=h^{0}(O(K-\alpha ))\neq 0$\textit{, that is either }$%
\alpha $\textit{\ or }$K-\alpha $\textit{\ is effective. Hence, there must
exist some integers }$m$\textit{'s such that }$\alpha +mK$\ \textit{is
effective because }$-K$\textit{\ is effective, we denote the largest such }$%
m $ \textit{as }$m_{\alpha }$\textit{. }

\textit{We claim that for every }$\alpha \in \Delta $\textit{, }$C_{\alpha
}\in |\alpha +m_{\alpha }K|$ \textit{is an irreducible }$(-2)$\textit{%
-curve. If so, then }$C=\cup _{\alpha \in \Delta }C_{\alpha }$\textit{\ is a
maximal }$ADE$\textit{\ curve of type }$\mathfrak{g}$. \textit{If there
exists reducible }$C_{\alpha }$\textit{, we write }$C_{\alpha }=\cup D_{i}$. 
\textit{Then each }$D_{i}$ \textit{is perpendicular to }$K$ as $-K$ \textit{%
is nef and }$C_{\alpha }\cdot K=0$\textit{. Since }$C_{\alpha }+K$ \textit{%
is not effective, every }$D_{i}+K$ \textit{is also not effective and }$%
D_{i}\notin |-K|$\textit{. Hence }$D_{i}^{2}=-2$\textit{\ for any }$i$%
\textit{\ as }$D_{i}^{2}=0$ \textit{will imply }$D_{i}\in |-K|$\textit{. We
know }$C_{\alpha }$\textit{\ is connected, this is because if }$C_{\alpha }$ 
\textit{is not connected, then one of its connected component must have
self-intersection zero from }$C_{\alpha }^{2}=-2$\textit{, which contradicts
to }$C_{\alpha }+K$\textit{\ is not effective. Hence }$C=\cup _{\alpha \in
\Delta }C_{\alpha }$ \textit{is an }$($\textit{affine}$)$ $ADE$ \textit{%
curve by Lemma \ref{curve2}. It is obvious that this curve strictly contains
a }$\mathfrak{g}$-\textit{curve, which contradicts to }$\mathcal{E}_{0}^{%
\widehat{E}_{8}}$ \textit{being maximal }$\mathfrak{g}$\textit{-deformable}.%
\textit{\ }

$\mathit{(iii)}$\textit{\ The proof is similar to }$\mathit{(ii)}$\textit{.}
\end{proof}

\begin{remark}
If $X_{9}$ admits an elliptic fibration, then we can find $m$ such that $%
h^{1}(X_{9},O(-mK))\neq 0$. Conversely, if $h^{1}(X_{9},O(-mK))\neq 0$, we
need to add the condition of $-K$ being nef to show that $X$ admits an
elliptic fibration. To see this, we take $x_{1},\cdots ,x_{5}$ to be $5$
points on a line $l\subset \mathbb{P}^{2}$, and another $4$ generic points $%
( $not on $l)$ $x_{6},\cdots ,x_{9}$ in $\mathbb{P}^{2}$. Then we have an
one parameter family of conics $C_{t}$'s passing through these $4$ points.
If we blow up $\mathbb{P}^{2}$ at these $9$ points and denote the strict
transforms of $l$ and $C_{t}$ with same notations, then $l^{2}=-4$, $%
C_{t}^{2}=0$. Moreover $C_{t}+l\in |-K|$ and $h^{0}(X_{9},O(-K))=2$. But $-K$
is not nef as $(-K)\cdot l=-2$, which implies that $X_{9}$ is not elliptic.
\end{remark}

From the above, we can easily deduce similar results for the $E_{n}$-bundle $%
\mathcal{E}_{0}^{E_{n}}$ over $X_{n}$ when $n\leq 8$, namely

(i) $\mathcal{E}_{0}^{E_{n}}$ is totally non-deformable if and only if the $%
n $ blowup points in $\mathbb{P}^{2}$ are in general position.

(ii) When $-K_{X_{n}}$ nef, $\mathcal{E}_{0}^{E_{n}}$ is maximal $\mathfrak{g%
}$-deformable if and only if $X_{n}$ has a maximal $\mathfrak{g}$ curve.

\section{Appendix}

In this appendix, we recall some results on affine Lie algebras \cite{Kac}%
\cite{LXZ}. If $(\mathfrak{g,}$ $[,])$ is a finite dimensional simple Lie
algebra, then the corresponding loop Lie algebra is $L\mathfrak{g:=g\otimes }%
\mathbb{C[}t,t^{-1}]$, with the Lie bracket defined by $[a\otimes
t^{n},b\otimes t^{m}]_{L\mathfrak{g}}=[a,b]\otimes t^{m+n}$, where $a,b\in 
\mathfrak{g}$, $m,n\in \mathbb{Z}$.

The corresponding untwisted affine Lie algebra $\widehat{\mathfrak{g}}$ is
constructed as a central extension of $L\mathfrak{g}$, with one-dimensional
center $\mathbb{C}c$, i.e. $\widehat{\mathfrak{g}}=L\mathfrak{g\oplus }%
\mathbb{C}c$. The Lie bracket on $\widehat{\mathfrak{g}}$ is defined by the
formula $[a\otimes t^{n}+\lambda c,b\otimes t^{m}+\mu c]_{L\mathfrak{g}%
}=[a,b]\otimes t^{m+n}+n\delta _{n+m,0}k(a,$ $b)c$, where $\lambda ,\mu \in 
\mathbb{C}$ and $k$ is the Killing form on $\mathfrak{g}$.

We can obtain the affine Dynkin diagram of $\widehat{\mathfrak{g}}$ from the
Dynkin diagram of $\mathfrak{g}$ by adding one node to it, corresponding to
the extended root and labelling as $C_{0}$. But in the affine $ADE$ except
affine $E_{8}$ case, from the symmetry of the affine Dynkin diagrams, we
have different choices of labelling the extended root.

The Institute of Mathematical Sciences and Department of Mathematics, The
Chinese University of Hong Kong, Shatin, N.T., Hong Kong

E-mail address: yxchen@math.cuhk.edu.hk

E-mail address: leung@math.cuhk.edu.hk

\end{document}